\pdfoutput=1
\documentclass[a4paper, 12pt]{amsart}
\usepackage[utf8]{inputenc}
\usepackage[margin=0.8in]{geometry}
\usepackage{amsmath,amssymb,enumerate,enumitem,nccmath}
\usepackage{bbm} %indicator function sembolünü kullanmak için
\usepackage{float}
\usepackage{graphics}
\usepackage{cite} %cite ederken sıralı ve [1-4] gibi yazması için.
\usepackage[table]{xcolor}
\usepackage{tikz}
\usetikzlibrary{patterns,arrows.meta,calc} 
\usepackage{adjustbox}
\usepackage{varwidth}
\usepackage{tcolorbox}
\tcbuselibrary{skins,breakable}
\usepackage{booktabs,array}
\renewcommand\arraystretch{1.25}
\usepackage{subfigure}
\usepackage{colortbl}
\usepackage{graphicx}
\usepackage{enumitem}
\usepackage[unicode, colorlinks, pdfpagelabels] {hyperref}
\usepackage{pgfplots}
\pgfplotsset{compat=1.17}
\usetikzlibrary{arrows}
\usetikzlibrary{decorations.text}
\usetikzlibrary{decorations.markings}
\usepackage{multirow}
\usepackage{pbox}
\usepackage{etoolbox}
\pgfplotsset{%
  myaddplot/.style = {%
    mark=*,
    mark size=.3pt,
    thick
  }
}
\usepackage{tabularx}
\definecolor{pastelgreen}{HTML}{C2D69B}
\definecolor{color1}{rgb}{0.55, 0.71, 0.0}
\definecolor{color2}{rgb}{0.63, 0.36, 0.94}
\definecolor{color3}{rgb}{0.01, 0.31, 0.59}

\newcommand{\mm}{\mathfrak{m}_{\mathtt{m}}}
\newcommand{\vfi}{\vec{\Phi}}
\newcommand{\scx}{\sin cx}

\newcommand{\smx}{\sin mx}
\newcommand{\vecb}{\vec{e}_{c,1}}

\newcommand{\ucb}{u_{c,1}}
\newcommand{\valphcb}{\vec{\alpha}_{c,1}}
\newcommand{\alpcbb}{\alpha_{c,1,1}}
\newcommand{\alpcbi}{\alpha_{c,1,2}}

\newcommand{\abs}[1]{\left\lvert #1 \right\rvert}
\newcommand{\norm}[1]{\left\lVert #1 \right\rVert}

\renewcommand\arraystretch{1.5}
\renewcommand\tabcolsep{3pt}

\newcommand{\indic}{\scalebox{1.25}{$\mathbbm{1}$}}
\newcommand{\ci}{\mathrm{i}}

\newtheorem{theorem}{Theorem}[section]

\newtheorem{remark}{Remark}[section]

\newtheorem{definition}{Definition}[section]

\definecolor{lime}{HTML}{A6CE39}
\DeclareRobustCommand{\orcidicon}{%
	\begin{tikzpicture}
	\draw[lime, fill=lime] (0,0)
	circle [radius=0.16]
	node[white] {{\fontfamily{qag}\selectfont \tiny ID}};
	\draw[white, fill=white] (-0.0625,0.095)
	circle [radius=0.007];
	\end{tikzpicture}
	\hspace{-2mm}
}

\foreach \x in {A, ..., Z}{%
	\expandafter\xdef\csname orcid\x\endcsname{\noexpand\href{https://orcid.org/\csname orcidauthor\x\endcsname}{\noexpand\orcidicon}}
}

\makeatletter
\g@addto@macro{\endabstract}{\@setabstract}
\newcommand{\authorfootnotes}{\renewcommand\thefootnote{\@fnsymbol\c@footnote}}%
\makeatother

\title[TRANSITION DYNAMICS IN GENERAL TWO-COMPONENT REACTION-DIFFUSION SYSTEMS]{Transition Dynamics in General Two-Component Reaction-Diffusion Systems}

\pgfplotstableread{%
d slope
1.020000 -48.03847806
1.040400 -23.78371026
1.061208 -15.69977767
1.082432 -11.65855555
1.104081 -9.23441883
1.126162 -7.61882575
1.148686 -6.46525834
1.171659 -5.60045755
1.195093 -4.92816830
1.218994 -4.39063749
1.243374 -3.95111312
1.268242 -3.58509381
1.293607 -3.27561704
1.319479 -3.01056667
1.345868 -2.78105754
1.372786 -2.58042572
1.400241 -2.40357526
1.428246 -2.24654262
1.456811 -2.10619865
1.485947 -1.98004004
1.515666 -1.86604026
1.545980 -1.76254123
1.576899 -1.66817322
1.608437 -1.58179477
1.640606 -1.50244704
1.673418 -1.42931869
1.706886 -1.36171857
1.741024 -1.29905425
1.775845 -1.24081499
1.811362 -1.18655817
1.847589 -1.13589825
1.884541 -1.08849786
1.922231 -1.04406055
1.960676 -1.00232471
1.999890 -0.96305863
2.039887 -0.92605633
2.080685 -0.89113410
2.122299 -0.85812756
2.164745 -0.82688917
2.208040 -0.79728615
2.252200 -0.76919867
2.297244 -0.74251828
2.343189 -0.71714665
2.390053 -0.69299437
2.437854 -0.66997994
2.486611 -0.64802900
2.536344 -0.62707346
2.587070 -0.60705094
2.638812 -0.58790411
2.691588 -0.56958026
2.745420 -0.55203078
2.800328 -0.53521081
2.856335 -0.51907886
2.913461 -0.50359652
2.971731 -0.48872818
3.031165 -0.47444074
3.213035 -0.43548534
3.405817 -0.40061918
3.610166 -0.36928082
3.826776 -0.34100615
4.056383 -0.31540814
4.299766 -0.29216145
4.557752 -0.27099063
4.831217 -0.25166083
5.121090 -0.23397066
5.428355 -0.21774637
5.754057 -0.20283726
6.099300 -0.18911202
6.465258 -0.17645566
6.853174 -0.16476709
7.264364 -0.15395705
7.700226 -0.14394645
8.162239 -0.13466496
8.651974 -0.12604982
9.171092 -0.11804488
9.721358 -0.11059973
10.304639 -0.10366899
10.922918 -0.09721174
11.578293 -0.09119092
12.272990 -0.08557298
13.009370 -0.08032739
13.789932 -0.07542638
14.617328 -0.07084460
15.494367 -0.06655891
16.424029 -0.06254808
17.409471 -0.05879268
18.454039 -0.05527485
19.561282 -0.05197816
20.734959 -0.04888748
21.979056 -0.04598887
23.297799 -0.04326942
24.695667 -0.04071725
26.177407 -0.03832131
27.748052 -0.03607139
29.412935 -0.03395803
31.177711 -0.03197243
33.048374 -0.03010642
35.031276 -0.02835240
37.133153 -0.02670331
39.361142 -0.02515255
41.722811 -0.02369399
44.226179 -0.02232192
46.879750 -0.02103098
49.692535 -0.01981620
52.674087 -0.01867290
55.834532 -0.01759675
59.184604 -0.01658365
62.735680 -0.01562981
66.499821 -0.01473165
70.489810 -0.01388583
74.719199 -0.01308922
79.202351 -0.01233888
83.954492 -0.01163207
88.991762 -0.01096620
94.331267 -0.01033884
99.991143 -0.00974774
}\slopedata
\pgfplotstableread{%
lam a
1.0000924 0.005
1.0003695 0.01
1.0008315 0.015
1.0014786 0.02
1.0023109 0.025
1.003329 0.03
1.0045331 0.035
1.0059238 0.04
1.0075016 0.045
1.0092672 0.05
1.0112213 0.055
1.0133648 0.06
1.0156986 0.065
1.0182238 0.07
1.0209413 0.075
1.0238526 0.08
1.0269589 0.085
1.0302616 0.09
1.0337624 0.095
1.0374631 0.1
1.0413654 0.105
1.0454714 0.11
1.0497833 0.115
1.0543034 0.12
1.0590343 0.125
1.0639788 0.13
1.0691398 0.135
1.0745206 0.14
1.0801246 0.145
1.0859557 0.15
1.0920179 0.155
1.0983157 0.16
1.1048539 0.165
1.1116378 0.17
1.1186732 0.175
1.1259664 0.18
1.1335241 0.185
1.1413541 0.19
1.1494645 0.195
1.1578647 0.2
1.1665647 0.205
1.175576 0.21
1.1849111 0.215
1.1945842 0.22
1.2046114 0.225
1.2150106 0.23
1.2258026 0.235
1.2370108 0.24
1.2486628 0.245
1.2607902 0.25
1.2734306 0.255
1.2866286 0.26
1.3004378 0.265
1.3149239 0.27
1.3301687 0.275
1.3462759 0.28
1.3633809 0.285
1.3816658 0.29
1.4013858 0.295
1.4229186 0.3
}\branchIstable
\pgfplotstableread{%
lam a
1.0000924 0.005
1.0003694 0.01
1.0008312 0.015
1.0014778 0.02
1.002309 0.025
1.003325 0.03
1.0045257 0.035
1.0059111 0.04
1.0074812 0.045
1.0092361 0.05
1.0111756 0.055
1.0132999 0.06
1.015609 0.065
1.0181027 0.07
1.0207812 0.075
1.0236443 0.08
1.0266922 0.085
1.0299249 0.09
1.0333422 0.095
1.0369443 0.1
1.0407311 0.105
1.0447026 0.11
1.0488588 0.115
1.0531998 0.12
1.0577254 0.125
1.0624358 0.13
1.067331 0.135
1.0724108 0.14
1.0776754 0.145
1.0831246 0.15
1.0887586 0.155
1.0945774 0.16
1.1005808 0.165
1.106769 0.17
1.1131419 0.175
1.1196995 0.18
1.1264418 0.185
1.1333689 0.19
1.1404806 0.195
1.1477771 0.2
1.1552584 0.205
1.1629243 0.21
1.170775 0.215
1.1788103 0.22
1.1870304 0.225
1.1954353 0.23
1.2040248 0.235
1.2127991 0.24
1.2217581 0.245
1.2309018 0.25
1.2402302 0.255
1.2497434 0.26
1.2594412 0.265
1.2693238 0.27
1.2793912 0.275
1.2896432 0.28
1.30008 0.285
1.3107014 0.29
1.3215076 0.295
1.3324986 0.3
}\branchIpred
\pgfplotstableread{%
lam a
0.99981669 0.005
0.99926697 0.01
0.99835154 0.015
0.99707152 0.02
0.99542846 0.025
0.99342437 0.03
0.99106162 0.035
0.98834301 0.04
0.98527167 0.045
0.98185108 0.05
0.97808506 0.055
0.97397767 0.06
0.96953326 0.065
0.96475639 0.07
0.95965181 0.075
0.95422445 0.08
0.94847937 0.085
0.94242174 0.09
0.93605681 0.095
0.92938988 0.1
0.92242629 0.105
0.91517138 0.11
0.90763048 0.115
0.89980888 0.12
0.89171182 0.125
0.88334449 0.13
0.87471197 0.135
0.86581927 0.14
0.85667128 0.145
0.84727278 0.15
0.83762844 0.155
0.82774278 0.16
0.8176202 0.165
0.80726497 0.17
0.79668119 0.175
0.78587285 0.18
0.77484377 0.185
0.76359764 0.19
0.75213801 0.195
0.74046825 0.2
0.72859164 0.205
0.71651128 0.21
0.70423014 0.215
0.69175107 0.22
0.67907675 0.225
0.66620977 0.23
0.65315256 0.235
0.63990744 0.24
0.62647659 0.245
0.61286211 0.25
0.59906594 0.255
0.58508994 0.26
0.57093583 0.265
0.55660525 0.27
0.54209973 0.275
0.5274207 0.28
0.51256948 0.285
0.49754731 0.29
0.48235535 0.295
0.46699465 0.3
}\branchIIunstable
\pgfplotstableread{%
lam a
0.99981667 0.005
0.99926667 0.01
0.99835001 0.015
0.99706669 0.02
0.9954167 0.025
0.99340004 0.03
0.99101672 0.035
0.98826674 0.04
0.9851501 0.045
0.98166679 0.05
0.97781681 0.055
0.97360017 0.06
0.96901687 0.065
0.9640669 0.07
0.95875027 0.075
0.95306697 0.08
0.94701701 0.085
0.94060039 0.09
0.9338171 0.095
0.92666714 0.1
0.91915052 0.105
0.91126724 0.11
0.9030173 0.115
0.89440069 0.12
0.88541741 0.125
0.87606747 0.13
0.86635087 0.135
0.8562676 0.14
0.84581767 0.145
0.83500107 0.15
0.82381781 0.155
0.81226789 0.16
0.8003513 0.165
0.78806804 0.17
0.77541812 0.175
0.76240154 0.18
0.7490183 0.185
0.73526839 0.19
0.72115181 0.195
0.70666857 0.2
0.69181867 0.205
0.6766021 0.21
0.66101887 0.215
0.64506897 0.22
0.62875241 0.225
0.61206919 0.23
0.5950193 0.235
0.57760274 0.24
0.55981952 0.245
0.54166964 0.25
0.5231531 0.255
0.50426989 0.26
0.48502001 0.265
0.46540347 0.27
0.44542027 0.275
0.4250704 0.28
0.40435387 0.285
0.38327067 0.29
0.36182081 0.295
0.34000428 0.3
}\branchIIpred
\pgfplotstableread{%
lam a
1.2540124 -0.3
1.249789 -0.295
1.2455652 -0.29
1.2413411 -0.285
1.2371167 -0.28
1.232892 -0.275
1.2286669 -0.27
1.2244415 -0.265
1.2202158 -0.26
1.2159897 -0.255
1.2117632 -0.25
1.2075365 -0.245
1.2033094 -0.24
1.1990819 -0.235
1.1948541 -0.23
1.190626 -0.225
1.1863975 -0.22
1.1821687 -0.215
1.1779396 -0.21
1.1737101 -0.205
1.1694802 -0.2
1.1652501 -0.195
1.1610195 -0.19
1.1567887 -0.185
1.1525574 -0.18
1.1483259 -0.175
1.144094 -0.17
1.1398617 -0.165
1.1356291 -0.16
1.1313962 -0.155
1.1271629 -0.15
1.1229292 -0.145
1.1186953 -0.14
1.1144609 -0.135
1.1102262 -0.13
1.1059912 -0.125
1.1017558 -0.12
1.0975201 -0.115
1.093284 -0.11
1.0890475 -0.105
1.0848107 -0.1
1.0805736 -0.095
1.0763361 -0.09
1.0720982 -0.085
1.06786 -0.08
1.0636214 -0.075
1.0593825 -0.07
1.0551432 -0.065
1.0509036 -0.06
1.0466636 -0.055
1.0424233 -0.05
1.0381826 -0.045
1.0339415 -0.04
1.0297001 -0.035
1.0254583 -0.03
1.0212161 -0.025
1.0169736 -0.02
1.0127308 -0.015
1.0084875 -0.01
1.004244 -0.005
}\branchIIIstable
\pgfplotstableread{%
lam a
0.99575569 0.005
0.99151101 0.01
0.98726597 0.015
0.98302056 0.02
0.97877479 0.025
0.97452866 0.03
0.97028216 0.035
0.96603529 0.04
0.96178806 0.045
0.95754046 0.05
0.95329249 0.055
0.94904415 0.06
0.94479544 0.065
0.94054636 0.07
0.93629692 0.075
0.9320471 0.08
0.92779691 0.085
0.92354635 0.09
0.91929542 0.095
0.91504412 0.1
0.91079244 0.105
0.90654039 0.11
0.90228797 0.115
0.89803517 0.12
0.893782 0.125
0.88952845 0.13
0.88527452 0.135
0.88102022 0.14
0.87676554 0.145
0.87251048 0.15
0.86825505 0.155
0.86399924 0.16
0.85974304 0.165
0.85548647 0.17
0.85122952 0.175
0.84697218 0.18
0.84271447 0.185
0.83845637 0.19
0.83419789 0.195
0.82993903 0.2
0.82567978 0.205
0.82142015 0.21
0.81716014 0.215
0.81289974 0.22
0.80863895 0.225
0.80437778 0.23
0.80011623 0.235
0.79585428 0.24
0.79159195 0.245
0.78732923 0.25
0.78306612 0.255
0.77880262 0.26
0.77453873 0.265
0.77027445 0.27
0.76600978 0.275
0.76174472 0.28
0.75747926 0.285
0.75321342 0.29
0.74894718 0.295
0.74468055 0.3
}\branchIIIunstable
\pgfplotstableread{%
lam a
1.2546479 -0.3
1.2504038 -0.295
1.2461596 -0.29
1.2419155 -0.285
1.2376714 -0.28
1.2334272 -0.275
1.2291831 -0.27
1.224939 -0.265
1.2206949 -0.26
1.2164507 -0.255
1.2122066 -0.25
1.2079625 -0.245
1.2037183 -0.24
1.1994742 -0.235
1.1952301 -0.23
1.1909859 -0.225
1.1867418 -0.22
1.1824977 -0.215
1.1782535 -0.21
1.1740094 -0.205
1.1697653 -0.2
1.1655211 -0.195
1.161277 -0.19
1.1570329 -0.185
1.1527887 -0.18
1.1485446 -0.175
1.1443005 -0.17
1.1400563 -0.165
1.1358122 -0.16
1.1315681 -0.155
1.127324 -0.15
1.1230798 -0.145
1.1188357 -0.14
1.1145916 -0.135
1.1103474 -0.13
1.1061033 -0.125
1.1018592 -0.12
1.097615 -0.115
1.0933709 -0.11
1.0891268 -0.105
1.0848826 -0.1
1.0806385 -0.095
1.0763944 -0.09
1.0721502 -0.085
1.0679061 -0.08
1.063662 -0.075
1.0594178 -0.07
1.0551737 -0.065
1.0509296 -0.06
1.0466854 -0.055
1.0424413 -0.05
1.0381972 -0.045
1.0339531 -0.04
1.0297089 -0.035
1.0254648 -0.03
1.0212207 -0.025
1.0169765 -0.02
1.0127324 -0.015
1.0084883 -0.01
1.0042441 -0.005
0.99575587 0.005
0.99151174 0.01
0.9872676 0.015
0.98302347 0.02
0.97877934 0.025
0.97453521 0.03
0.97029108 0.035
0.96604695 0.04
0.96180281 0.045
0.95755868 0.05
0.95331455 0.055
0.94907042 0.06
0.94482629 0.065
0.94058215 0.07
0.93633802 0.075
0.93209389 0.08
0.92784976 0.085
0.92360563 0.09
0.9193615 0.095
0.91511736 0.1
0.91087323 0.105
0.9066291 0.11
0.90238497 0.115
0.89814084 0.12
0.8938967 0.125
0.88965257 0.13
0.88540844 0.135
0.88116431 0.14
0.87692018 0.145
0.87267605 0.15
0.86843191 0.155
0.86418778 0.16
0.85994365 0.165
0.85569952 0.17
0.85145539 0.175
0.84721125 0.18
0.84296712 0.185
0.83872299 0.19
0.83447886 0.195
0.83023473 0.2
0.8259906 0.205
0.82174646 0.21
0.81750233 0.215
0.8132582 0.22
0.80901407 0.225
0.80476994 0.23
0.8005258 0.235
0.79628167 0.24
0.79203754 0.245
0.78779341 0.25
0.78354928 0.255
0.77930515 0.26
0.77506101 0.265
0.77081688 0.27
0.76657275 0.275
0.76232862 0.28
0.75808449 0.285
0.75384035 0.29
0.74959622 0.295
0.74535209 0.3
}\branchIIIpred
\pgfplotstableread{%
t a
0 0.07
1.332 0.067548667
2.664 0.065745805
3.996 0.0640218
5.328 0.062251003
6.66 0.060409964
7.992 0.058499921
9.324 0.056528814
10.656 0.054506948
11.988 0.052445775
13.32 0.050357357
14.652 0.048254
15.984 0.046147914
17.316 0.044050911
18.648 0.04197414
19.98 0.03992787
21.312 0.037921322
22.644 0.035962562
23.976 0.034058444
25.308 0.032214603
26.64 0.030435486
27.972 0.028724414
29.304 0.027083672
30.636 0.025514609
31.968 0.024017751
33.3 0.022592908
34.632 0.021239289
35.964 0.019955602
37.296 0.018740154
38.628 0.017590934
39.96 0.016505688
41.292 0.015481992
42.624 0.014517303
43.956 0.01360901
45.288 0.012754472
46.62 0.011951054
47.952 0.011196148
49.284 0.010487198
50.616 0.0098217156
51.948 0.0091972902
53.28 0.0086115988
54.612 0.0080624116
55.944 0.0075475959
57.276 0.0070651176
58.608 0.006613042
59.94 0.0061895323
61.272 0.0057928487
62.604 0.0054213452
63.936 0.0050734674
65.268 0.0047477489
66.6 0.0044428075
67.932 0.0041573423
69.264 0.003890129
70.596 0.0036400171
71.928 0.0034059254
73.26 0.0031868386
74.592 0.0029818041
75.924 0.0027899278
77.256 0.0026103712
78.588 0.0024423484
79.92 0.0022851225
}\timeIIbelow
\pgfplotstableread{%
t a
0 0.1
1.332 0.099830308
2.664 0.10204244
3.996 0.10533147
5.328 0.10959873
6.66 0.11512101
7.992 0.12245366
9.324 0.13260012
10.656 0.14753399
11.988 0.17180634
13.32 0.2192732
14.652 0.37089091
}\timeIIabove
\pgfplotstableread{%
t a
0 -0.02
3.332 -0.022256198
6.664 -0.024606033
9.996 -0.027021151
13.328 -0.029469777
16.66 -0.031918327
19.992 -0.034333221
23.324 -0.036682703
26.656 -0.03893845
29.988 -0.041076812
33.32 -0.043079583
36.652 -0.044934269
39.984 -0.046633925
43.316 -0.048176622
46.648 -0.049564698
49.98 -0.050803887
53.312 -0.05190243
56.644 -0.052870256
59.976 -0.053718259
63.308 -0.054457716
66.64 -0.055099825
69.972 -0.055655377
73.304 -0.05613453
76.636 -0.056546668
79.968 -0.056900337
83.3 -0.057203224
86.632 -0.057462175
89.964 -0.057683239
93.296 -0.057871722
96.628 -0.058032255
99.96 -0.058168858
103.292 -0.058285007
106.624 -0.058383701
109.956 -0.058467516
113.288 -0.05853866
116.62 -0.058599026
119.952 -0.058650228
123.284 -0.058693645
126.616 -0.058730452
129.948 -0.058761648
133.28 -0.058788084
136.612 -0.058810483
139.944 -0.058829459
143.276 -0.058845534
146.608 -0.058859149
149.94 -0.058870681
153.272 -0.058880447
156.604 -0.058888717
159.936 -0.05889572
163.268 -0.058901651
166.6 -0.058906672
169.932 -0.058910924
173.264 -0.058914524
176.596 -0.058917572
179.928 -0.058920153
183.26 -0.058922337
186.592 -0.058924187
189.924 -0.058925753
193.256 -0.058927079
196.588 -0.058928202
199.92 -0.058929152
}\timeIIIleft
\pgfplotstableread{%
t a
0 0.02
1 0.021397999
2 0.022922121
3 0.024588048
4 0.026414143
5 0.028422035
6 0.030637383
7 0.033090855
8 0.035819425
9 0.038868074
10 0.042292092
11 0.046160207
12 0.050558942
13 0.055598771
14 0.061423034
15 0.068221149
16 0.076248751
17 0.085859418
18 0.097556565
19 0.11208226
20 0.13057773
21 0.15489355
22 0.18824249
23 0.23673645
24 0.31361231
25 0.45390177
26 0.79085196
27 2.7224549
}\timeIIIright

\begin{document}
\maketitle
\begin{center}

  \normalsize
  \authorfootnotes
Taylan Şengül\footnote{\href{taylan.sengul@marmara.edu.tr}{taylan.sengul@marmara.edu.tr}}\textsuperscript{1} \orcidA{}, Burhan Tiryakioglu\footnote{\href{burhan.tiryakioglu@marmara.edu.tr}{burhan.tiryakioglu@marmara.edu.tr}}\textsuperscript{1} \orcidB{} and Esmanur Yıldız Akıl\footnote{\href{esmanur.yildiz@yeditepe.edu.tr}{esmanur.yildiz@yeditepe.edu.tr} (Corresponding Author) }\textsuperscript{1,2} \orcidC{}
 \par \bigskip

  \textsuperscript{1}Department of Mathematics, Marmara University, 34722 Istanbul, Turkey \par
  \textsuperscript{2}Department of Mathematics, Yeditepe University, 34755 Istanbul, Turkey\par \bigskip

\today
\end{center}
\begin{abstract}
This work provides a comprehensive characterization of the first dynamic transition of a general two-component reaction–diffusion system on a bounded interval in one spatial dimension.
By rigorously addressing the difficulties arising from the non-linear interactions between the two components, this study significantly extends the scope of previous studies on dynamic transitions of reaction–diffusion systems, which focused on single-equation cases.
The analysis is conducted by employing center manifold reduction within the standard classification of dynamic transition theory.
A key contribution of this study is that the nonlinear terms are treated to involve arbitrary powers of the unknowns $u$ and $v$.
Furthermore, the theoretical framework is complemented by the analysis of a representative and simplified system.
As a result of this generalization, the dynamics of a broad class of fundamental reaction–diffusion systems can be systematically classified.
\end{abstract}

%\tableofcontents

\section{Introduction}\label{sec:introduction}

Reaction-diffusion systems are coupled nonlinear partial differential equation systems that play a fundamental role in modeling physical, chemical, biological processes \cite{bezekci2025spectral, cherniha2017nonlinear, volpert2025reaction}.
Reaction–diffusion systems have a wide variety of applications, including epidemiological models \cite{zhu2020spread}, spatially extended systems \cite{kolinichenko2024effect}, and nonlinear Burgers-type equations \cite{hao2025steady}.
A general two-component reaction-diffusion system defined in one spatial dimension is expressed as
\begin{equation}\label{rd-eqn-tr}
\frac{d \vec{u}}{dt} = L_{\lambda} \vec{u} + \vec{g}(\vec{u}).
\end{equation}
Here, $\vec{u}=(u,v)^{T}$ is the unknown vector, $L_{\lambda}$ is the linear diffusion operator, and $\vec{g}(\vec{u})$ represents the nonlinear interactions between the components.
In this framework, $\vec{u}$ can model various interactions such as chemical concentrations \cite{duzs2025reactor}, predator–prey relationships \cite{DU2025101621, YANG2026525}, or competition between two species \cite{cherniha2025reaction}.
Diffusion terms have different physical interpretations depending on the application.
For instance, in chemical systems they represent the tendency of molecules to randomly spread and homogenize, while in ecological models they demonstrate the spatial movement of species.
In general, diffusion mathematically describes the spatial spreading of matter or species in a system.

The interaction between nonlinear reaction terms and diffusion often leads the system to exhibit complex dynamical behavior.
For this reason, investigating the exchange of stabilities of reaction–diffusion equations and understanding the dynamic transitions that arise from changes in parameter values has been an intriguing research topic.
Dynamic transition theory, developed by Ma and Wang \cite{ptd}, has been a fundamental tool for analyzing transitions between various stable states of such systems and for classifying the types of these transitions.

Dynamic transition theory has been successfully applied to various models.
For example, dynamic transitions in three-dimensional reaction–diffusion systems \cite{li2021dynamic}, stability changes in competition models \cite{liu2022dynamic}, dynamic transition analysis for activator-substrate systems \cite{li2023dynamic}, and recent analyses of the general Brusselator model \cite{choi2024turing} are some of the significant examples considered within the scope of this theory.

The main goal of this study is to identify the dynamic transitions of the general two-variable reaction–diffusion system in one spatial dimension at the critical parameter values and to characterize the transition types.
This work extends our previous dynamic transition analyses for dissipative equations in one spatial dimension to a more general system consisting of coupled equations \cite{csengul2022dynamic, csengul2023dynamic, csengul2024first}.

In this paper the linear operator $L_{\lambda}$ is assumed to be self-adjoint.
For $L_{\lambda} = D \partial_{xx} + A$ with $D = \mathrm{diag}(d,1)$ and a constant reaction matrix $A$, this means that $A$ is symmetric.
Off-diagonal coupling is allowed, but the activator--inhibitor sign structure $a_{12}a_{21}<0$ of Turing systems is not.
We restrict to this case because the eigenvectors are then orthogonal, so the center manifold coefficients and the transition numbers can be computed in closed form.
A symmetric reaction matrix arises as the linearization of gradient kinetics $\vec{g} = \nabla F$, which includes multi-component Allen--Cahn \cite{fazly2025qualitative} and real Ginzburg--Landau \cite{malomed2022new} type systems.
Classical models such as the Brusselator \cite{zhao2025new}, Gray--Scott \cite{zheng2018existence, gandy2022analyzing, choi2025turing, al2021unified}, Schnakenberg \cite{yang2022cross, alfifi2022stability}, Turing \cite{sun2021turing, chen2025pattern} and Gierer--Meinhardt \cite{song2017pattern} systems therefore lie outside the present hypotheses.
We plan to address the non-self-adjoint case, together with these applications, in a forthcoming paper; see \autoref{sec:conclusions}.

The main assumptions of this work are the following.
The linear operator $L_{\lambda}$ is assumed to be general differential operator with real eigenvalues and pure sine eigenmodes.
A typical example is the Laplacian operator with Dirichlet boundary conditions.
In particular, the critical eigenmode is assumed to be of the form $\vec{e}_{c,1} = \vec{\alpha}_{c,1} \sin cx$ where $\vec{\alpha}_{c,1} \in \mathbb{R}^2$ is a constant vector.
We assume a simple real critical eigenvalue $\beta_{c,1}$ with wavenumber $c \in \mathbb{Z}_{\geq 1}$ which becomes unstable at some $\lambda = \lambda_0$.

The nonlinear operator is assumed to be a higher order polynomial of the unknown variables $u$ and $v$ with real coefficients.
This is usually the case in applications after taking the perturbation around the equilibrium solution.

Under these assumptions, we obtain, under a generic nondegeneracy condition, an exact formulation for the system's first transition number $\mathcal{P}$ and the critical index $\mathfrak{m} \in \mathbb{Z}_{\geq 2}$ in terms of the coefficients of the nonlinear operator and the eigenvalues of the linear operator which completely describe the system's first dynamic transition.
In summary, if $\mathfrak{m}$ is even the first transition is mixed type and if $\mathfrak{m}$ is odd the first transition is continuous if $\mathcal{P} <0$ and jump type if $\mathcal{P} >0$.

We further show that the first transition number $\mathcal{P}$ in general consists of two parts: a self-interaction term and a cross-interaction term.
The self-interaction term is due to the nonlinear interaction of the critical mode $\vecb$ with itself, while the cross-interaction term is due to the nonlinear interaction of the critical mode $\vecb$ with the center manifold $\vec{\Phi}$.
\emph{One of the main contributions of this work is the exact formulation of $\mathcal{P}$ and $\mathfrak{m}$ in terms of the coefficients of the nonlinear operator and the eigenvalues of the linear operator.}

We also give an application of our main results, first in the case where the reaction matrix is a multiple of the identity matrix and then for a one-parameter family with off-diagonal coupling, and show the possible transition scenarios depending on the system data.

The main result is obtained in three steps.
By the PES condition a single mode $\vec{e}_{c,1}$ loses stability at $\lambda_0$, and the center manifold reduction gives the scalar equation \eqref{reduced-eqn} for its amplitude, whose leading term determines $\mathcal{P}$ and $\mathfrak{m}$.
Projecting the nonlinearity on the center manifold onto the critical mode splits this coefficient into the self and cross interaction numbers.
Both are computed in closed form for every order $N$, which is what allows nonlinearities of arbitrary order.

This paper is organized as follows. In \autoref{sec:formulation-of-the-problem}, we present the general framework
of the problem together with the assumptions used in this work. In \autoref{sec:main-results} we introduce
key definitions such as the transition number and the critical index. The main theorem of the paper is also
stated in this section. \autoref{sec:application} provides an application of the theoretical results to reaction–diffusion systems.
The proof of the main theorem is given in \autoref{sec:proof}. Finally, we sum up the results and discuss future work
in \autoref{sec:conclusions}.

\section{Formulation of the Problem}\label{sec:formulation-of-the-problem}

We consider the following system of two partial differential equations of the form
\begin{equation}\label{main-eqn}
    \begin{bmatrix}
        \frac{\partial u}{\partial t} \\
        \frac{\partial v}{\partial t}
    \end{bmatrix} = L_{\lambda}
    \begin{bmatrix}
        u \\
        v
    \end{bmatrix} + \vec{g}(u,v),
\end{equation}
with unknown functions $u(x,t)$ and $v(x,t)$ on the spatial domain
\[
x \in (0,\pi) \subset \mathbb{R},
\]
$t \ge 0$ denotes the time and $\lambda$ is a real control parameter.

Now we will give the assumptions on the linear operator $L_{\lambda}$ and the nonlinear operator $\vec{g}(u,v)$.

\subsection{The assumptions on the linear operator $L_{\lambda}$}
We assume that the linear operator $L_{\lambda}$ is self-adjoint so that its spectrum consist only of real eigenvalues
\[
\{ \beta_{m,n} \in \mathbb{R}: m \in \mathbb{Z}_{\geq 1}, n \in \{1,2\} \},
\]
with corresponding eigenfunctions given by
\begin{equation*}
    \vec{e}_{m,n}= \vec{\alpha}_{m,n} \sin mx,
\end{equation*}
where
\begin{equation}\label{eigenfunctions}
    \vec{\alpha}_{m,n}=  \begin{bmatrix}
        \alpha_{m,n,1} \\ \alpha_{m,n,2}
    \end{bmatrix}, \quad m=1,2, \ldots  \quad \text{and} \quad n=1,2.
\end{equation}
We note that the eigenfunctions in this study consist of pure sine modes.
As a typical example, this is the case when the linear operator is $L_{\lambda} = \frac{d^2}{dx^2} + \lambda $ with homogeneous Dirichlet boundary conditions.
\begin{equation}\label{boundary-conditions}
    u = v = 0 \quad \text{at } x = 0 \text{ and } x = \pi.
\end{equation}
For example, in this case, the linear operator is defined on $H^2(0,\pi) \cap H_0^1(0,\pi)$.
However, we allow the linear operator to be of higher order and not necessarily second order.

We assume that the vectors $\vec{\alpha}_{m,n}$ are normalized so that
\begin{equation}\label{normalization}
    \norm{\vec{\alpha}_{m, n}}=1 \quad \text{for all} \quad m\ge 1, \quad n =1,2,
\end{equation}
and, since $L_{\lambda}$ is self-adjoint, $\vec{\alpha}_{m,1} \perp \vec{\alpha}_{m,2}$, so that the eigenfunctions satisfy the orthogonality condition
\[
\langle \vec{e}_{m,n}, \vec{e}_{m',n'} \rangle = \frac{\pi}{2} \langle \vec{\alpha}_{m,n}, \vec{\alpha}_{m',n'} \rangle =  \frac{\pi}{2} \delta_{m,m'} \delta_{n,n'}, \qquad \forall m,m' \in \mathbb{Z}_{\geq 1}, \quad n,n' \in \{1,2\},
\]
where $\delta$ stands for the Kronecker delta.

We assume that the linear operator $L_{\lambda}$ satisfies the principle of exchange of stabilities (PES) condition
\begin{equation}\label{PES-condition}
    \begin{aligned}
    & \beta_{c,1} \begin{cases}
        <0, \quad \lambda < \lambda_0 \\
        =0, \quad \lambda = \lambda_0 \\
        >0, \quad \lambda > \lambda_0 \\
    \end{cases} \\
    &\beta_{m,n}(\lambda_0)<0, \qquad \forall (m,n) \neq (c,1),
    \end{aligned}
\end{equation}
for some critical value $\lambda_0 \in \mathbb{R}$ of the control parameter and for some \textbf{critical wavenumber} $c \in \mathbb{Z}_{\ge 1}$.
Throughout, $L_{\lambda}$ is a sectorial operator with compact resolvent, so that $\beta_{m,n} \to -\infty$ as $m \to \infty$.
Together with \eqref{PES-condition}, this gives $\sup_{(m,n) \ne (c,1)} \beta_{m,n}(\lambda_0) < 0$, which is the spectral gap required by the center manifold theorem \cite[Theorem~A.1.1]{ptd}.
The above PES condition means that largest (the critical) eigenvalue $\beta_{c,1}$ is simple with corresponding critical mode
\[
\vec{e}_{c,1}= \vec{\alpha}_{c,1} \scx \text{ for some } c \in \mathbb{Z}_{\geq 1}.
\]
For second order linear operators, $c$ is usually $1$, but for higher order linear operators, $c$ can be greater than $1$.
See for example~\cite{csengul2022dynamic}.

An example for the linear operator satisfying the above assumptions is given in \autoref{sec:application}.

\subsection{The assumptions on the nonlinear operator $\vec{g}(u,v)$}

The nonlinear operator $\vec{g}(u,v)$ is given by
\begin{equation}\label{nonlinear-operator}
    \vec{g}(u,v) = \sum_{\substack{k+l=2 \\ k, l \in \mathbb{N}_0}}^{\infty} \vec{a}_{k,l} u^k v^l,
\end{equation}
where
\begin{equation*}
\vec{a}_{k,l}=\begin{bmatrix}
    a_{k,l,1} \\ a_{k,l,2}
  \end{bmatrix}.
\end{equation*}
That is we consider arbitrary powers of the unknown variables $u$ and $v$ with arbitrary real coefficients $\vec{a}_{k,l}$.

We assume that $\limsup_{k+l \to \infty} \norm{\vec{a}_{k,l}}^{1/(k+l)} < \infty$, that is, the power series \eqref{nonlinear-operator} has a positive radius of convergence.
Since $H_0^1(0,\pi)$ is a Banach algebra which embeds continuously in $C[0,\pi]$, the series converges in $H_0^1(0,\pi)$ near the origin.
Hence $\vec{g} : H_1 \to H$ is analytic, where $H = L^2(0,\pi)^2$ and $H_1 \subset H_0^1(0,\pi)^2$ is the domain of $L_{\lambda}$.

\subsection{Background: the reduced system and the type of transition}
The PES condition implies that the system undergoes a dynamic transition at $\lambda = \lambda_0$ from the trivial solution $\vec{u}_c \equiv 0$.
A crucial part is then to determine the type of transition that occurs at this critical value. This is the main goal of this work.

To determine the type of transition, it is necessary to study the solution on the center manifold by writing it as
\[
\vec{u}(x,t) = u_{c,1}(t) \vec{e}_{c,1}(x) + \vec{\Phi}(x,t),
\]
where $u_{c,1}(t)$ is the amplitude of the critical mode and $\Phi(x,t)$ is the center manifold spanned by the stable modes.

By the center manifold reduction, it is well known that the first transition of the system is governed by the truncated \textbf{reduced equation}
\begin{equation}\label{reduced-eqn}
    \frac{d u_{c,1}}{dt} = \beta_{c,1} u_{c,1} + \mathcal{P} u_{c,1}^{\mathfrak{m}},
\end{equation}
where $\mathcal{P} \ne 0$ is the \textbf{transition number} and $\mathfrak{m} \in \mathbb{Z}_{\geq 2}$ is the \textbf{critical index}.

By the standard classification of the dynamic transition theory \cite[Chapter~2, Theorems~2.3.1 and~2.3.2]{ptd}, the following theorem holds.
\begin{theorem}\label{theorem mP}
Assume that all the above assumptions are satisfied.
Then the system \eqref{main-eqn} undergoes a dynamic transition at $\lambda = \lambda_0$ from the trivial solution $\vec{u}_c \equiv 0$.
The type of transition is determined as follows:
\begin{enumerate}
    \item[i)] If $\mathfrak{m}$ is odd and $\mathcal{P} < 0$ then there is a Type-I (continuous) transition.
    \item[ii)] If $\mathfrak{m}$ is odd and $\mathcal{P} > 0$ then there is a Type-II (jump) transition.
    \item[iii)] If $\mathfrak{m}$ is even then there is a Type-III (mixed) transition.
\end{enumerate}
\end{theorem}
Here, Type-I transition means that the bifurcating states are near the trivial solution, Type-II transition means that the system leaves a local neighborhood of the trivial solution, and Type-III transition means that a local neighborhood of the trivial solution can be written as a disjoint union of two open sets, where one of them the system moves to a local steady state and the other one the system moves out of the neighborhood.

\subsection{The main goal of this work}
The main goal of this work is to determine the transition number $\mathcal{P}$ and the critical index $\mathfrak{m}$ in terms of the coefficients of the nonlinear operator $\vec{a}_{k,l}$, the eigenvalues $\beta_{m,n}$ and the eigenvectors $\vec{\alpha}_{m,n}$ of the linear operator.

We show that the transition number $\mathcal{P}$ may consist of two different contributions: the self interaction number and the cross interaction number.
The \textbf{self interaction number} consists of the nonlinear interactions of the critical mode with itself.
The \textbf{cross interaction number} consists of the nonlinear interactions of the critical mode with other (stable) modes.

In our previous study \cite{csengul2024first}, we only considered the first transition of a single equation.
This is the first result where the first transition is characterized in a system of more than one equation.
This considerably complicates the analysis of the transition number $\mathcal{P}$ and the critical index $\mathfrak{m}$.

Finally, we provide an example of a reaction diffusion system to illustrate the main results.

\section{Main Results}\label{sec:main-results}

For the following definitions, we let $\indic_{A}$ denote the indicator function of a set $A$ defined as
\begin{equation}\label{indicator-func}
    \indic_{A}(x)=
    \begin{cases}
        1, &\quad x \in A \\
        0, &\quad x \notin A.
    \end{cases}
\end{equation}

We also assume that $0^0=1$ whenever it appears in the expressions.

\begin{definition}\label{S-definition}
\textbf{Nth self interaction number} is defined by
    \begin{equation}\label{S_N-last}
        \mathcal{S}_{N}= \mathit{s}_{N,c} \sum_{\substack{k+l=N \\ 0 \leq k, l \leq N}}^{} (\vec{a}_{k,l} \cdot \vec{\alpha}_{c,1}) \alpha_{c,1,1}^k \alpha_{c,1,2}^l,
    \end{equation}
where
\begin{equation}\label{s_Nc}
    \mathit{s}_{N,c} = \frac{2}{\pi} B\left(\frac{N+2}{2},\frac{1}{2}\right) \left[\indic_{\mathit{2\mathbb{Z}+1}} (N) + \frac{1}{c} \indic_{\mathit{2\mathbb{Z}}} (N) \indic_{\mathit{2\mathbb{Z}+1}} (c) \right].
\end{equation}
Here $B$ is the so called beta-function.
\end{definition}

The number $\mathcal{S}_N$ measures the interaction of the critical mode with itself: it is the term of order $N$ in the projection of $\vec{g}(u_{c,1}\vec{e}_{c,1})$ onto $\vec{e}_{c,1}$, and $\mathit{s}_{N,c} = \frac{2}{\pi}\int_0^{\pi}\sin^{N+1}cx\,dx$.
For even $N$ this integral vanishes unless $c$ is odd, which is the parity rule expressed by the indicator terms in \eqref{s_Nc}.

\begin{definition}\label{C-definition}
\textbf{Nth cross interaction number} is defined for $N \ge 2$ by
\begin{equation}\label{C_N-last}
	\mathcal{C}_{N}= \sum_{\delta \in \Lambda_N}^{} \frac{1}{\beta_{m,n}} c_{\delta} (\vec{a}_{k,l} \cdot \vec{\alpha}_{c,1}) (\vec{a}_{k',l'} \cdot \vec{\alpha}_{m,n}),
\end{equation}
where
\begin{equation}\label{index-set}
\begin{aligned}
    \Lambda_{N} = \big\{ \,
        &\delta = (N_1, N_2, k, l, k', l', m, n) \in \mathbb{Z}_{\geq 0}^{8} :\
        N_1 \ge 2,\, N_2 \ge 2,\, N_1+N_2=N+1, \\
        & k+l=N_1,\, k'+l'=N_2,\, m \geq 1,\, n \in \{1,2\},\, (m,n) \neq (c,1)
    \big\},
\end{aligned}
\end{equation}
is the index set and
\begin{equation*}
    c_{\delta}= -\frac{4}{\pi^2} \left( k \alpha^{k+k'-1}_{c,1,1} \alpha^{l+l'}_{c,1,2} \alpha_{m,n,1} + l \alpha^{k+k'}_{c,1,1} \alpha^{l+l'-1}_{c,1,2} \alpha_{m,n,2} \right) \mathcal{J}_{m,c}^{1,N_1} \mathcal{J}_{m,c}^{1,N_2} , \quad \delta  \in \Lambda_{N},
\end{equation*}
with
\begin{align}\label{Int-2}
    \mathcal{J}_{m,c}^{1,N} = & \frac{1}{(2 \ci)^{N+1}} \sum_{j=0}^{N} (-1)^j \binom{N}{j} \left[  \indic_{\{0\}}(\gamma_j)  \pi
    - \indic_{\{2m\}}(\gamma_j) \pi + \indic_{\mathbb{Z} \setminus {\{0, 2m\} }}(\gamma_j) \dfrac{2m \left( 1-(-1)^{\gamma_{j}} \right)}{\ci \gamma_{j} (\gamma_{j} - 2m)} \right],
   \end{align}
and
\begin{equation*}
    \gamma_j=m+(N-2j)c.
\end{equation*}
\end{definition}

The number $\mathcal{C}_N$ measures an indirect interaction: through the nonlinear terms of order $N_2$ the critical mode excites the stable mode $(m,n)$, whose amplitude on the center manifold carries the factor $1/\beta_{m,n}$, and this mode acts back on the critical mode through the terms of order $N_1$.
The integrals $\mathcal{J}^{1,N}_{m,c}=\int_0^{\pi}\sin mx\,\sin^{N}cx\,dx$ are the corresponding spatial overlaps, and \eqref{Int-2} gives them in closed form for every $N$.
Since $N_1,N_2 \ge 2$, the index set $\Lambda_2$ is empty and $\mathcal{C}_2=0$.

We define the \textbf{Nth transition number} by
\begin{equation} \label{PN}
\mathcal{P}_N = \mathcal{S}_{N} + \mathcal{C}_{N}, \qquad N \ge 2.
\end{equation}

Before stating the main theorem we record the smallest cases of the definitions, which are the ones used in the application of \autoref{sec:application}.

\begin{remark}[The smallest cases]\label{rem:small-cases}
Evaluating the beta-function in \eqref{s_Nc} gives $B(2,\tfrac12)=\tfrac43$ and $B(\tfrac52,\tfrac12)=\tfrac{3\pi}{8}$, so that
    \[
        \mathit{s}_{2,c}=\frac{8}{3\pi c}\,\indic_{2\mathbb{Z}+1}(c), \qquad \mathit{s}_{3,c}=\frac34 ,
    \]
the latter independent of $c$.
Consequently
    \[
        \mathcal{S}_2 = \frac{8}{3\pi c}\,\indic_{2\mathbb{Z}+1}(c) \sum_{k+l=2} (\vec{a}_{k,l} \cdot \vec{\alpha}_{c,1}) \alpha_{c,1,1}^{k} \alpha_{c,1,2}^{l}, \qquad
        \mathcal{S}_3 = \frac34 \sum_{k+l=3} (\vec{a}_{k,l} \cdot \vec{\alpha}_{c,1}) \alpha_{c,1,1}^{k} \alpha_{c,1,2}^{l} .
    \]
For the cross term $N_1=N_2=2$ is the only splitting, so that
    \begin{align*}
        \mathcal{C}_3 = -\frac{4}{\pi^{2}} \sum_{(m,n) \ne (c,1)} \frac{\left[\mathcal{J}_{m,c}^{1,2}\right]^{2}}{\beta_{m,n}}
        \sum_{k+l=2} \sum_{k'+l'=2} & (\vec{a}_{k,l} \cdot \vec{\alpha}_{c,1}) (\vec{a}_{k',l'} \cdot \vec{\alpha}_{m,n}) \\
        & \times \left( k \alpha^{k+k'-1}_{c,1,1} \alpha^{l+l'}_{c,1,2} \alpha_{m,n,1} + l \alpha^{k+k'}_{c,1,1} \alpha^{l+l'-1}_{c,1,2} \alpha_{m,n,2} \right) .
    \end{align*}
When in addition $c=1$ the indicator in $\mathit{s}_{2,c}$ is satisfied and $\mathit{s}_{2,1}=8/(3\pi)$; this is the situation of \autoref{sec:application}.
\end{remark}

Let $\mm$ denote the lowest degree $k+l$ with $\vec{a}_{k,l} \ne \vec{0}$.
As in \cite{csengul2024first}, we assume a \textbf{genericity condition}
\begin{equation}\label{genericity}
    \mathcal{P}_N \ne 0 \quad \text{for some } 2 \le N \le 2\mm-1 ,
\end{equation}
which fails only on the zero set of a nonzero polynomial in the coefficients $\vec{a}_{k,l}$; in \autoref{example-main-theorem} it fails exactly on the line $\mathcal{P}_3=0$ of \autoref{fig:P_figure1}.

Under the assumptions stated in the previous section we are ready to state the main theorem of this work.
\begin{theorem}\label{main-theorem}
Suppose that \eqref{genericity} holds.
Then the \textbf{critical index} $\mathfrak{m}$ in \eqref{reduced-eqn} is given by
\[
\mathfrak{m} = \min \{N \ge 2: \mathcal{P}_N \neq 0\},
\]
and the \textbf{transition number} in \eqref{reduced-eqn} is given by
\[
\mathcal{P} = \mathcal{P}_{\mathfrak{m}}.
\]
\end{theorem}

The proof of the above theorem, given in \autoref{sec:proof}, depends on the idea that the reduced equation \eqref{reduced-eqn} can be written as
\begin{equation*}
        \frac{d u_{c,1}}{dt} =
        \beta_{c,1} u_{c,1}
        + \sum_{N=2}^{2\mm-1} (\mathcal{S}_{N} + \mathcal{C}_{N}) u_{c,1}^{N} + O\big(u_{c,1}^{2\mm}\big).
\end{equation*}

Some remarks are now in order.
\begin{remark}~
    \begin{enumerate}
    \item When $c$ is even, $\mathcal{P}_{2N}=0$ and the transition is Type–I or Type–II.
    Indeed, the reflection $x\mapsto\pi-x$ commutes with the system and maps $\sin cx$ to $-\sin cx$, so the reduced equation is odd in $u_{c,1}$; under \eqref{genericity}, $\mathfrak{m}$ is odd and the result follows from \autoref{theorem mP}.
    \item
    Now suppose that either $c$ is odd or the nonlinear operator has a term of the form $\vec{a}_{k,l} \ne \vec{0}$ for some $k+l$ odd. If we also suppose that the critical vector $\vec{\alpha}_{c,1}$ has both components $\alpha_{c,1,1}$ and $\alpha_{c,1,2}$ non-zero, then the self interaction part $\mathcal{S}$ is generically non-zero.
\end{enumerate}
\end{remark}

\section{Proof of \autoref{main-theorem}}\label{sec:proof}
In this section, we will prove \autoref{main-theorem} by analyzing the dynamics around $\lambda= \lambda_0$ when the solutions are close to trivial solution $\vec{u} \equiv 0$.

The steps were outlined in \autoref{sec:introduction}.
The center manifold reduction gives a scalar equation for the amplitude $u_{c,1}$.
We split its nonlinear part $g_{c,1}$ into the terms free of the center manifold function $\vec{\Phi}$, which give $\mathcal{S}$, and the terms linear in $\vec{\Phi}$, which give $\mathcal{C}$.
Since $\vec{\Phi} = O(u_{c,1}^{\mm})$, the discarded terms, quadratic in $\vec{\Phi}$ or of higher order in the center manifold approximation, are $O(u_{c,1}^{2\mm})$ and do not affect $\mathcal{P}_N$ for $N \le 2\mm-1$.
For this purpose, we represent the solution on the center manifold as
\begin{equation}\label{cm-eqn}
    \begin{aligned}
        \vec{u}= u_{c,1}(t) \vec{e}_{c,1} + \vfi,
    \end{aligned}
\end{equation}
where the center manifold function is
\begin{equation}\label{cm-funct}
    \begin{aligned}
        &\vfi= \begin{bmatrix}
            \Phi_{u} \\ \Phi_{v}
        \end{bmatrix}= \sum_{(m,n) \neq (c,1)}^{} u_{m,n}(t) \vec{e}_{m,n}(x).
    \end{aligned}
\end{equation}
The above sum is over all indices $(m,n) \neq (c,1)$ where $m \in \mathbb{Z}_{\geq 1}$ and $n \in \{1,2\}$.
Thus the center manifold is spanned by all the stable modes except the critical mode $\vec{e}_{c,1}$.

Equivalently to equations \eqref{cm-eqn} and \eqref{cm-funct}, we can express the solution as
\begin{equation*}
    \begin{aligned}
        &u= u_{c,1}(t) \alpcbb \scx + \Phi_u \\
        &v= u_{c,1}(t) \alpcbi \scx + \Phi_v,
    \end{aligned}
\end{equation*}
where
\begin{equation*}
    \begin{aligned}
        &\Phi_u=\sum_{(m,n) \neq (c,1)}^{} u_{m,n}(t) \alpha_{m,n,1} \smx \\
        &\Phi_v=\sum_{(m,n) \neq (c,1)}^{} u_{m,n}(t) \alpha_{m,n,2} \smx.
    \end{aligned}
\end{equation*}

The center manifold function satisfies the tangency condition $\Phi= O(u_{c,1}^2) $, $u_{c,1} \to 0$.
However, depending on the nonlinearity, the tangency can be of higher order as we will make clear in the following proof.

Center manifold analysis near the critical point $(\vec{u}, \lambda)= (0, \lambda_0)$ provides a means to obtain the scalar reduced equation given by
\begin{equation*}
    \frac{d u_{c,1}}{dt} = f(u_{c,1}).
\end{equation*}
The reduced equation captures the dynamics of the full PDE system when $\abs{\lambda - \lambda_0}$ and $\abs{u_{c,1}}$ are both small.
Here
\begin{equation*}
    f(u_{c,1}) = \frac{1}{\langle \vecb, \vecb \rangle }\langle L_{\lambda} \vec{u} + \vec{g}(\vec{u}), \vecb \rangle = \beta_{c,1} u_{c,1} + g_{c,1},
  \end{equation*}
where
  \begin{equation*}
    g_{c,1}= \frac{1}{\langle \vecb, \vecb \rangle } \int_{0}^{\pi} \vec{g}(u_{c,1}(t) \vec{e}_{c,1} + \vfi) \cdot \vecb \, dx = O(\abs{ u_{c,1} }^{2}),
\end{equation*}
as given in \cite{csengul2024first}.
Since
\begin{equation*}
    \vecb= \valphcb \scx,
\end{equation*}
$g_{c,1}$ can be shown as
\begin{equation}\label{gc1}
    g_{c,1}= \frac{2}{\pi} \frac{1}{\norm{\valphcb}^{2}} \int_{0}^{\pi} \vec{g} (\ucb(t)\vecb + \vfi) \cdot \valphcb   \scx \, dx.
\end{equation}

If we rearrange the equation \eqref{gc1}, it has the form
\begin{equation*}
    g_{c,1}= \frac{2}{\pi} \frac{1}{\norm{\valphcb}^{2}} \int_{0}^{\pi} \sum_{k+l=2}^{\infty} \vec{a}_{k,l} \cdot \valphcb (\ucb \alpcbb \scx + \Phi_u)^{k}(\ucb \alpcbi \scx + \Phi_v)^{l}  \scx \, dx.
\end{equation*}

By expanding the powers in the above integrand to at most linear order in $\vec{\Phi}$, and using the normality $\norm{\valphcb}^{2} = 1$, see \eqref{normalization}, we obtain
\begin{align}\label{gm1-last-eqn}
    g_{c,1} &= \frac{2}{\pi} \sum_{k+l=2}^{\infty} \vec{a}_{k,l} \cdot \valphcb \Bigg[ \int_{0}^{\pi} \ucb^{k+l} (t) \alpcbb^{k} \alpcbi^{l} \sin^{k+l+1} cx \, \, dx \nonumber  \\
    &\quad + \int_{0}^{\pi} \ucb^{k+l-1}(t)\alpcbb^{k-1} \alpcbi^{l-1} \sin^{k+l} cx \left( l \alpcbb \Phi_v + k \alpcbi \Phi_u \right) \, dx \Bigg]  + \text{h.o.t.}
\end{align}
where $\text{h.o.t.}$ denotes higher order terms in $\vec{\Phi}$.
We now want to express $g_{c,1}$ as
\begin{equation*}
    g_{c,1} = \mathcal{S} + \mathcal{C} + O(\vec{\Phi}^{2}).
\end{equation*}

Here $\mathcal{S}$ includes terms that are independent of the center manifold $\vec{\Phi}$ and $\mathcal{C} $ consists of linear terms in $\vec{\Phi}$.
Thus $\mathcal{S}$ is formed by the self-nonlinear interaction of the critical mode $\vecb$ and $\mathcal{C}$ is formed by the nonlinear interactions of the critical mode $\vecb$ with the center manifold $\vec{\Phi}$.

We introduce the notation
\begin{equation}\label{int-defn}
    \mathcal{J}_{k_1 , k_2, \ldots , k_n }^{l_1, l_2, \ldots , l_n} = \int_{0}^{\pi}\sin ^{l_1} (k_1 x)  \sin (k_2 x) ^{l_2} \ldots \sin (k_n x) ^{l_n}  \, dx.
\end{equation}

\subsection{Computation of self interaction term $\mathcal{S}$}
From \eqref{gm1-last-eqn}, the self interaction term $\mathcal{S}$ is given by
\begin{equation}\label{S-value}
    \mathcal{S}= \frac{2}{\pi} \sum_{k+l=2}^{\infty} \vec{a}_{k,l} \cdot \valphcb  \alpcbb^{k} \alpcbi^{l}  \mathcal{J}_{c}^{k+l+1} \ucb^{k+l} = \sum_{N = 2}^{\infty} S_N u_{c,1}^N,
\end{equation}
where
\begin{equation*}
    S_N = \frac{2}{\pi} \sum_{k+l=N}^{} \vec{a}_{k,l} \cdot \valphcb  \alpcbb^{k} \alpcbi^{l}  \mathcal{J}_{c}^{N+1}, \quad N \geq 2.
\end{equation*}
By utilizing \eqref{int-defn}, we express the integral as
\begin{equation*}
    \mathcal{J}_{c}^{N+1} = \int_{0}^{\pi} \sin^{N+1} cx \, dx \, = \frac{1}{c} \int_{0}^{c \pi} \sin^{N+1} \theta \, d \theta \, = \frac{1}{c} \sum_{n=0}^{c-1}\int_{n \pi}^{(n+1)\pi} \sin^{N+1} \theta \, d\theta.
\end{equation*}
By the periodicity and symmetry properties of sine function, we can express the integral as
\begin{equation*}
    \mathcal{J}_{c}^{N+1} = \frac{1}{c} \sum_{n=0}^{c-1} a_n \int_{0}^{\pi} \sin^{N+1} \theta \, d\theta,
\end{equation*}
where $a_n = 1$ if $n$ is even and $a_n = (-1)^{N+1}$ if $n$ is odd.
Considering the parity of $N$ and $c$ gives the following cases:
\begin{equation}\label{Jc-cases}
    \mathcal{J}_{c}^{N+1}=
    \begin{cases}
        0, &\text{$N$ is even and $c$ is even} \\[6pt]
        \frac{1}{c} \int_{0}^{\pi} \sin^{N+1} \theta \, d \theta, &\text{$N$ is even and $c$ is odd} \\[6pt]
        \int_{0}^{\pi} \sin^{N+1} \theta \, d \theta, &\text{$N$ is odd}.
    \end{cases}
\end{equation}
To evaluate $\int_{0}^{\pi} \sin^{N+1} \theta \, d \theta$, we split the integral as
\begin{equation}\label{int-0pi}
    \int_{0}^{\pi} \sin^{N+1} \theta \, d \theta = \int_{0}^{\frac{\pi}{2}} \sin^{N+1} \theta \, d \theta + \int_{\frac{\pi}{2}}^{\pi} \sin^{N+1} \theta \, d \theta \, = 2 \int_{0}^{\frac{\pi}{2}} \sin^{N+1} \theta \, d \theta.
\end{equation}
Next, we recall the identity that establishes the relationship between the Gamma and Beta functions, see \cite{abramowitz1964handbook},
\begin{equation}\label{beta-gamma-funct}
    B(a,b) = \frac{\Gamma(a)\Gamma(b)}{\Gamma(a+b)} = 2 \int_{0}^{\frac{\pi}{2}} \cos^{2a-1}\theta \sin^{2b-1}\theta \, d\theta.
\end{equation}
By choosing $a=\frac{1}{2}$, $b=\frac{N+2}{2}$ in \eqref{beta-gamma-funct},
\begin{equation*}
    B\left(\frac{1}{2}, \frac{N+2}{2}\right) = \sqrt{\pi} \frac{\Gamma\left(\frac{N+2}{2}\right)}{\Gamma\left(\frac{N+3}{2}\right)} = 2 \int_{0}^{\frac{\pi}{2}} \sin^{N+1}\theta \, d\theta.
\end{equation*}
Combining the above result with \eqref{Jc-cases}, the integral can be written in the closed form as
\begin{equation}\label{Int-2-Bform}
    \mathcal{J}_{c}^{N+1} = B\left(\frac{N+2}{2},\frac{1}{2}\right)
    \begin{cases}
    \frac{\left(c+1-(-1)^{N}(c-1)\right)}{2c} &\text{ if } c \equiv 1 \text{   } (\bmod 2) \bigskip\\
    \frac{\left(1-(-1)^{N}\right)}{2} &\text{ if } c \equiv 0 \text{   } (\bmod 2).
    \end{cases}
    \end{equation}
By making use of the indicator function, the above expression can be rewritten as
\begin{equation}\label{int-eqn1}
    \mathcal{J}_{c}^{N+1} = B\left(\frac{N+2}{2},\frac{1}{2}\right) \left[ \indic_{\mathit{2\mathbb{Z}+1}} (N) + \frac{1}{c} \indic_{\mathit{2\mathbb{Z}}} (N) \indic_{\mathit{2\mathbb{Z}+1}} (c) \right].
\end{equation}
Here $\mathit{2\mathbb{Z}}$ represents the set of even integers while $\mathit{2\mathbb{Z}+1}$ demonstrate the set of odd integers.

Finally by defining,
\[
    s_{N,c} = \frac{2}{\pi} \mathcal{J}_{c}^{N+1},
\]
we obtain the expression for $\mathcal{S}_N$ as in \autoref{S-definition}.

\subsection{Computation of cross interaction term $\mathcal{C}$}

Next, we compute the part of \eqref{gm1-last-eqn}, which contains linear terms in $\vec{\Phi}$ which is called the cross interaction term $\mathcal{C}$.
It is given as
\begin{align}\label{C-value}
    \mathcal{C}= \frac{2}{\pi} \sum_{N_1 = 2}^{\infty} \sum_{k+l=N_1} \vec{a}_{k,l} \cdot \valphcb  \alpcbb^{k-1} \alpcbi^{l-1} \ucb^{N_1-1} \int_{0}^{\pi} \left( l \alpcbb \Phi_v \sin^{N_1} cx + k \alpcbi \Phi_u \sin^{N_1} cx \right)\, dx.
\end{align}
To determine a formulation for $\mathcal{C}$, let us substitute $\Phi_u$ and $\Phi_v$ into \eqref{C-value}.
Using the integral notation \eqref{int-defn}, $\mathcal{C}$ takes the form of
\begin{align*}
    \mathcal{C}= \sum_{N_1=2}^{\infty} \hat{c}_{N_1} \ucb^{N_1-1},
\end{align*}
where
\begin{equation*}
    \hat{c}_{N_1} = \frac{2}{\pi}\sum_{k+l=N_1} \sum_{(m,n) \neq (c,1)}^{} \vec{a}_{k,l} \cdot \valphcb \alpcbb^{k-1} \alpcbi^{l-1}  u_{m,n} \mathcal{J}_{m, c}^{1, N_1} \left( k \alpcbi  \alpha_{m,n,1} + l \alpcbb \alpha_{m,n,2} \right).
\end{equation*}
To obtain a closed relation for $\mathcal{C}$, we need to compute the coefficients $u_{m,n}$ of the center manifold function $\vec{\Phi}$ corresponding to the stable mode $\vec{e}_{m,n}$.

For this we employ the center manifold approximation for critical crossing of real eigenvalues, see \cite[Theorem~A.1.1]{ptd},
\begin{equation*}
    \langle   -P_{\text{stable}} L_{\lambda} \vfi \, , \, \vec{e}_{m,n} \rangle = \langle \vec{g}(\vec{u}_{c,1}) \, , \, \vec{e}_{m,n} \rangle + \text{h.o.t.}, \quad \forall (m,n) \neq (c,1),
\end{equation*}
where $P_{\text{stable}}$ is the projection onto the stable space spanned by the modes except the critical mode.
By the orthogonality of the eigenfunctions, we have
\begin{align*}
    -\sum_{(m',n')\neq (c,1)} \int_{0}^{\pi} \beta_{m',n'} u_{m',n'}(t) \vec{e}_{m',n'} \cdot \vec{e}_{m,n} \, dx = - u_{m,n} \beta_{m,n} \frac{\pi}{2} = \int_{0}^{\pi} \vec{g}(\vec{u}_{c,1}) \cdot \vec{e}_{m,n} \, dx.
\end{align*}

Using the integral notation \eqref{int-defn} for the right hand side, $u_{m,n}$ finally takes the form
\begin{equation}\label{u-mn}
    u_{m,n}= -\frac{1}{\beta_{m,n}} \frac{2}{\pi} \sum_{N_2=2}^{\infty}\ucb^{N_2}  \sum_{k'+l'=N_2}^{\infty} \vec{a}_{k',l'} \cdot \vec{\alpha}_{m,n} \alpha^{k'}_{c,1,1} \alpha^{l'}_{c,1,2}  \mathcal{J}_{m,c}^{1,N_2}.
\end{equation}
Plugging \eqref{u-mn} into $\hat{c}_{N_1}$ and collecting the powers of $\ucb$, the term $\ucb^{N_1-1}\ucb^{N_2}$ contributing to the power $N=N_1+N_2-1$, we obtain
\begin{align*}
    \mathcal{C}= \sum_{N=3}^{\infty} \mathcal{C}_{N} \ucb^{N},
\end{align*}
where the Nth cross interaction number $\mathcal{C}_{N}$ is as in \autoref{C-definition}.
Now, we obtain an exact relation for the integral $\mathcal{J}_{m,c}^{1,N}$ which is given below.
\begin{equation*}
    \mathcal{J}_{m, c}^{1, N} = \int_{0}^{\pi} \sin mx \sin^{N} cx \, dx.
\end{equation*}
By using Euler's identity, we obtain
\begin{align}\label{int-calc}
    \mathcal{J}_{m,c}^{1,N} &= \frac{1}{(2 \ci)^{N+1}} \sum_{j=0}^{N} (-1)^j \binom{N}{j} \int_{0}^{\pi} (e^{\ci \gamma_{j} x} - e^{\ci (\gamma_j - 2m) x} )\, dx,
\end{align}
where
\begin{equation*}
    \gamma_j=m+(N-2j)c.
\end{equation*}

The last integral is evaluated as
\begin{equation}\label{int-J}
    \int_{0}^{\pi} (e^{\ci \gamma_{j} x} - e^{\ci (\gamma_j - 2m) x} )\, dx = \begin{cases}
        \pi,  \quad &\text{$\gamma_j = 0$} \\
        - \pi,  \quad &\text{$\gamma_j = 2m$} \\
        \dfrac{2m \left( 1-(-1)^{\gamma_{j}} \right)}{\ci \gamma_{j} (\gamma_{j} - 2m)},  \quad &\text{$\gamma_j \neq 0$ and $\gamma_j \neq 2m$}. \\
    \end{cases}
\end{equation}

By making use of indicator function, we obtain the expression for $\mathcal{J}_{m,c}^{1,N}$ as in \eqref{Int-2}.

\section{Application to Reaction Diffusion Systems}\label{sec:application}
To illustrate the main results, we consider the following reaction–diffusion system
\begin{equation}\label{ex-main-eqn}
    \frac{d}{dt}
    \begin{bmatrix}
        u \\ v
    \end{bmatrix}
    =
    \begin{bmatrix}
        d \,\partial_{xx} & 0 \\
        0 & \partial_{xx}
    \end{bmatrix}
    \begin{bmatrix}
        u \\ v
    \end{bmatrix}
    +
    A
    \begin{bmatrix}
        u \\ v
    \end{bmatrix}
    +
    \begin{bmatrix}
        g_{1}(u,v) \\ g_{2}(u,v)
    \end{bmatrix}.
\end{equation}
This framework covers reaction--diffusion systems with a symmetric reaction matrix $A$, that is, those whose linearization at the equilibrium is self-adjoint.
See \autoref{sec:conclusions} for the non-self-adjoint case.

The constant matrix $A$ represents the linear reaction dynamics, that is, the Jacobian of the reaction terms (excluding diffusion) evaluated at the equilibrium.
The nonlinear terms $g_1$ and $g_2$ which describe the nonlinear reaction dynamics are assumed to be of the form \eqref{nonlinear-operator}.
That is
\begin{equation}\label{ex-nonlinear-operator}
    \begin{aligned}
        g_1 &= a_{2,0,1} u^2 + a_{1,1,1} uv + a_{0,2,1} v^2 + \ldots,\\
        g_2 &= a_{2,0,2} u^2 + a_{1,1,2} uv + a_{0,2,2} v^2 + \ldots
    \end{aligned}
\end{equation}

We assume that the system is subject to homogeneous Dirichlet boundary conditions
\begin{equation}\label{ex-boundary-conditions}
    u(0)=u(\pi)=v(0)=v(\pi)=0.
\end{equation}

\subsection{Diagonal reaction matrix}\label{subsec:diagonal}
For simplicity of the discussion, we first assume that the linear reaction matrix $A$ is of the form
\begin{equation}\label{ex-linear-reaction-matrix}
A=\begin{bmatrix}\lambda & 0 \\[2pt] 0 & \lambda\end{bmatrix}.
\end{equation}

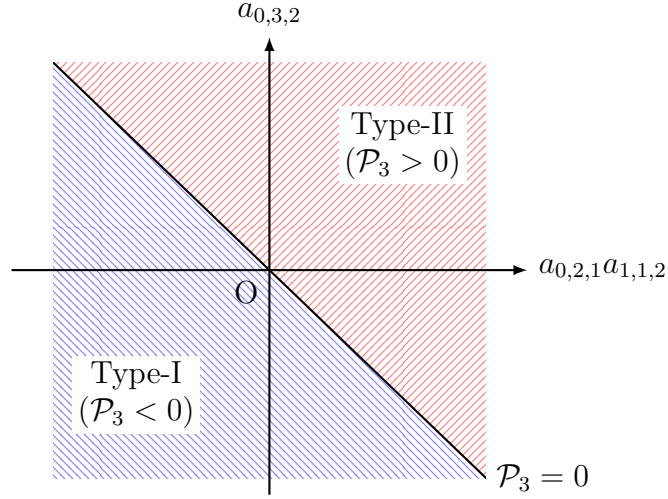
\begin{figure}[htbp]
    \centering
    \begin{tikzpicture}[>=latex,scale=1.1]
        % slope -0.96 corresponds to d=2 (Table 1)
        \begin{scope}
            \clip (-2.6,-2.5) rectangle (2.6,2.5);
            \fill[pattern=north east lines, pattern color=red!50]
                (-3,{-0.96*(-3)}) -- (3,{-0.96*3}) -- (3,3) -- (-3,3) -- cycle;
            \fill[pattern=north west lines, pattern color=blue!50]
                (-3,{-0.96*(-3)}) -- (3,{-0.96*3}) -- (3,-3) -- (-3,-3) -- cycle;
            \draw[thick] (-3,{-0.96*(-3)}) -- (3,{-0.96*3});
        \end{scope}
        \draw[->,thick] (-3.1,0)--(3.1,0) node[right] {$a_{0,2,1} a_{1,1,2}$};
        \draw[->,thick] (0,-2.7)--(0,2.8) node[above] {$a_{0,3,2}$};
        \node[right] at (2.6,{-0.96*2.6}) {$\mathcal{P}_3=0$};
        \node[fill=white, inner sep=2pt, align=center] at (1.6,1.5) {Type-II\\($\mathcal{P}_3>0$)};
        \node[fill=white, inner sep=2pt, align=center] at (-1.6,-1.5) {Type-I\\($\mathcal{P}_3<0$)};
        \node[below left] at (0,0) {O};
    \end{tikzpicture}
    \caption{The regions of the $(a_{0,2,1}a_{1,1,2},\,a_{0,3,2})$-plane where the transition number $\mathcal{P}_3$ is positive or negative, for $d=2$; the line $\mathcal{P}_3=0$ has slope $-\tfrac43 k(d)$.}
    \label{fig:P_figure1}
\end{figure}

The following theorem is a direct consequence of \autoref{main-theorem}.
\begin{theorem}\label{example-main-theorem}
Suppose that $d>1$.
Then the system \eqref{ex-main-eqn}-\eqref{ex-linear-reaction-matrix} exhibits a first dynamic transition at $\lambda = 1$ from the trivial solution.
Then the transition types up to $\mathfrak{m}=3$ can be classified as follows.
    \begin{enumerate}
    \item If $a_{0,2,2}\ne0$ then $\mathcal{P} = \mathcal{S}_2 = \frac{8}{3\pi}a_{0,2,2}$, $\mathfrak{m}=2$ and the transition is Type-III (Mixed).
    \item If $a_{0,2,2}=0$ then there exists a positive real number $k(d)$ which depends on the diffusion value $d$ such that $\mathcal{P}_3 = k(d) a_{0,2,1}a_{1,1,2} + \frac{3}{4} a_{0,3,2}$ and the transition scenarios are as shown in \autoref{fig:P_figure1}. In particular, we have the following:
    \begin{enumerate}
        \item If $ \mathcal{P}_3 > 0$ then $\mathcal{P} = \mathcal{P}_3$, $\mathfrak{m}=3$ and the transition is Type-II (Jump).
        \item If $ \mathcal{P}_3 < 0$ then $\mathcal{P} = \mathcal{P}_3$, $\mathfrak{m}=3$ and the transition is Type-I (Continuous).
        \item If $ \mathcal{P}_3 = 0$ then $\mathfrak{m}> 3$ and the coefficients of higher order nonlinear terms are needed to determine the transition type.
    \end{enumerate}
    \end{enumerate}
\end{theorem}

When $a_{0,2,2}=0$, $a_{0,3,2} \ne 0$, $a_{0,2,1}a_{1,1,2}\ne0$, generically $\mathcal{P}_3 \ne 0$ which means that $\mathcal{P}_3 = 0$ is a set of measure zero which is a line in the $(a_{0,2,1}a_{1,1,2},a_{0,3,2})$-plane as shown in \autoref{fig:P_figure1}.

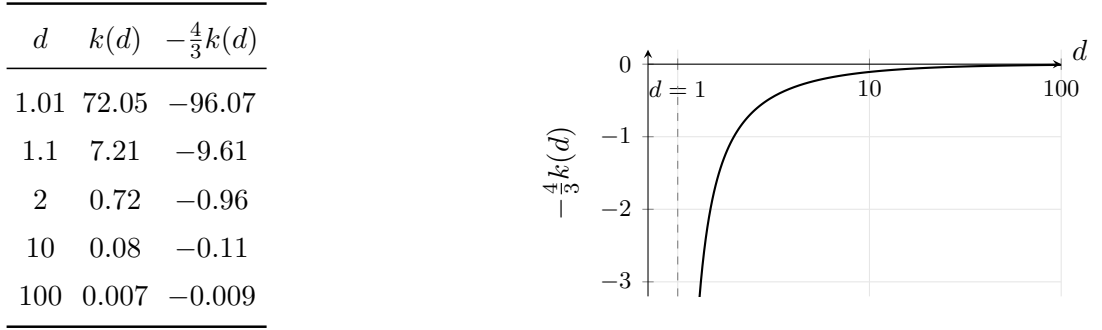
\begin{figure}[htbp]
    \centering
    \begin{minipage}[c]{0.42\textwidth}
        \centering
        \renewcommand{\arraystretch}{1.4}
        \small
        \begin{tabular}{ccc}
        \toprule
        $d$ & $k(d)$ & $-\tfrac43 k(d)$ \\
        \midrule
        1.01 & 72.05 & $-96.07$ \\
        1.1  & 7.21  & $-9.61$ \\
        2    & 0.72  & $-0.96$ \\
        10   & 0.08  & $-0.11$ \\
        100  & 0.007 & $-0.009$ \\
        \bottomrule
        \end{tabular}
    \end{minipage}\hfill
    \begin{minipage}[c]{0.52\textwidth}
        \centering
        \begin{tikzpicture}
            \begin{axis}[
                width=0.8\linewidth, height=0.55\linewidth,
                xlabel={$d$}, ylabel={$-\tfrac43 k(d)$},
                xmode=log, log basis x=10,
                xmin=0.7, xmax=100, ymin=-3.2, ymax=0.2,
                axis x line=middle, axis y line=left,
                xlabel style={at={(axis description cs:1,1)}, anchor=west},
                ylabel style={font=\small},
                xtick={1,10,100}, xticklabels={$d=1$,$10$,$100$}, xticklabel pos=upper,
                grid=major, grid style={gray!20},
                tick label style={font=\scriptsize}, label style={font=\small},
            ]
                \addplot[thick, black] table[x=d, y=slope] {\slopedata};
                \addplot[dashed, gray] coordinates {(1,-3.2) (1,-0.3)};
            \end{axis}
        \end{tikzpicture}
    \end{minipage}
    \caption{The constant $k(d)$ and the slope $-\tfrac43 k(d)$ of the line $\mathcal{P}_{3}=0$ in \autoref{fig:P_figure1}: computed values for selected $d$ (left) and the slope for $d>1$ on a logarithmic axis (right), where the dashed line marks the vertical asymptote at $d=1$.}
    \label{SlopeP-figure}
\end{figure}

Since $k(d) > 0$, the slope of the line $\mathcal{P}_3=0$ in \autoref{fig:P_figure1} is always negative.
In fact, the slope of the line $\mathcal{P}=0$ is given by $- \dfrac{4}{3} k(d)$ which we compute for some diffusion values $d>1$ in \autoref{SlopeP-figure}.
To compute the exact values of $k(d)$, we need to compute $\mathcal{J}_{m,1}^{1,2}$ for all $m \in \mathbb{Z}_{\geq 1}$.
A simplified formula for $\mathcal{J}_{m,1}^{1,2}$, derived in the proof below, is
\begin{equation}\label{J-simplified}
    \mathcal{J}_{m,1}^{1,2}=
    \begin{cases}
        \dfrac{4}{m(4-m^2)}, & m \text{ odd}, \\[6pt]
        0, & m \text{ even}.
    \end{cases}
\end{equation}
Formula \eqref{J-simplified} also explains the behavior of $k(d)$ for $d$ near $1$ and for large $d$.
As $d \to 1^{+}$, the $m=1$ term of \eqref{kd} dominates, so that $k(d) \sim \frac{4}{\pi^{2}}\,\frac{(4/3)^{2}}{d-1} \approx \frac{0.7205}{d-1}$.
As $d \to \infty$, we may replace $1/(1-dm^{2})$ by $-1/(dm^{2})$ in each term, since the series converges uniformly in $d$ as shown in the proof of \autoref{example-main-theorem}.
This gives $k(d) \sim C/d$ with $C = \frac{4}{\pi^{2}} \sum_{m \ge 1} \frac{[\mathcal{J}_{m,1}^{1,2}]^{2}}{m^{2}} \approx 0.7237$.
The two constants are nearly equal, so $k(d) \approx 0.72/(d-1)$ is a good approximation for all $d>1$.
Furthermore, \autoref{SlopeP-figure} illustrates how the slope of the line $\mathcal{P}_3=0$ varies with the diffusion value $d>1$.
This plot shows that the slope has a vertical asymptote at $d=1$ and approaches zero as $d$ becomes large.
This implies that the line $\mathcal{P}_3=0$ in \autoref{fig:P_figure1} becomes nearly horizontal for large diffusion values.
Finally, the results in \autoref{table:combined} confirm the consistency between the theoretical predictions and the computed parameter values.

\begin{table}[htbp]
    \centering
    \renewcommand{\arraystretch}{1.4}
    \small

    \begin{adjustbox}{max width=\linewidth}
    \begin{tabular}{c|cc|cc}
    \hline
    \textbf{Scenario} &
    \multicolumn{2}{c|}{\textit{A: Type-I (Continuous)}} &
    \multicolumn{2}{c}{\textit{B: Type-II (Jump)}} \\
    \midrule
    \textbf{Parameters} &
    \multicolumn{2}{c|}{$a_{0,2,1}=1$, $a_{1,1,2}=-2$, $a_{0,3,2}=-3$} &
    \multicolumn{2}{c}{$a_{0,2,1}=2$, $a_{1,1,2}=3$, $a_{0,3,2}=4$} \\
    \midrule
    \textbf{\boldmath$d$ Values} & \boldmath{$\mathcal{C}_3$} & \boldmath{$\mathcal{P}$} & \boldmath{$\mathcal{C}_3$} & \boldmath{$\mathcal{P}$} \\
    \midrule
    1.01 & $-144.11$ & $-146.36$ & $432.33$ & $435.33$ \\
    1.1  & $-14.42$  & $-16.67$  & $43.25$  & $46.25$ \\
    2    & $-1.44$   & $-3.69$   & $4.33$   & $7.33$ \\
    10   & $-0.16$   & $-2.41$   & $0.48$   & $3.48$ \\
    100  & $-0.01$   & $-2.26$   & $0.04$   & $3.04$ \\
    \bottomrule
    \end{tabular}
    \end{adjustbox}

    \vspace{4pt}
    \caption{Comparison of $\mathcal{C}_3$ and $\mathcal{P}$ values for Scenarios A and B with $\mathcal{S}_3 = -\dfrac{9}{4}$ and $\mathcal{S}_3 = 3$, respectively.}
    \label{table:combined}
    \end{table}

\subsubsection*{Proof of \autoref{example-main-theorem}}
We will use our main theorem (\autoref{main-theorem}) to prove \autoref{example-main-theorem}.

\subsubsection*{The Linear Problem} We first investigate the linear problem.
The linear operator is
\begin{equation*}
    L_\lambda =
    \begin{bmatrix}
        d \,\partial_{xx} & 0 \\
        0 & \partial_{xx}
    \end{bmatrix}
    + A,
\end{equation*}
with Dirichlet boundary conditions \eqref{ex-boundary-conditions}.

The associated eigenvalue problem becomes
\begin{equation*}
    \begin{aligned}
        &\beta\, u = d\, u_{xx} + \lambda u\\
        &\beta\, v = v_{xx} + \lambda v,
    \end{aligned}
\end{equation*}
with separated solutions
\begin{equation*}
    u = u_k \sin kx,
    \qquad
    v = v_k \sin kx.
\end{equation*}

For $d \neq 1$, the eigenvalues and corresponding eigenvector components for $k \geq 1$ are
\begin{equation*}
    \begin{aligned}
        &\beta_{k,1}= \lambda - k^2, \quad u_k = 0, \quad v_k = 1,\\
        &\beta_{k,2}= \lambda - d k^2, \quad u_k = 1, \quad v_k = 0.
    \end{aligned}
\end{equation*}
For $d>1$, the first few eigenvalues and eigenvectors of the linear operator take the form
\begin{equation*}
    \begin{aligned}
        \beta_{1,1} &= \lambda - 1, \quad
        e_{1,1} = \begin{bmatrix} 0 \\ 1 \end{bmatrix} \sin x, \quad
        \beta_{2,1} = \lambda - 4, \quad
        e_{2,1} = \begin{bmatrix} 0 \\ 1 \end{bmatrix} \sin 2x, \quad \ldots \\
        \beta_{1,2} &= \lambda - d, \quad
        e_{1,2} = \begin{bmatrix} 1 \\ 0 \end{bmatrix} \sin x, \quad
        \beta_{2,2} = \lambda - 4d, \quad
        e_{2,2} = \begin{bmatrix} 1 \\ 0 \end{bmatrix} \sin 2x, \quad \ldots
    \end{aligned}
\end{equation*}
Therefore, by \eqref{eigenfunctions}, for $m \in \mathbb{Z}_{\geq 1}$,
\begin{equation*}
    \alpha_{m,1}= \begin{bmatrix} 0 \\ 1 \end{bmatrix},
    \qquad
    \alpha_{m,2}= \begin{bmatrix} 1 \\ 0 \end{bmatrix}.
\end{equation*}

We verify the PES condition \eqref{PES-condition}.
Since $d>1$, we have $\beta_{m,2}=\lambda-dm^{2} < \lambda-m^{2}=\beta_{m,1}$ and $\beta_{m,1}$ decreases in $m$, so the largest eigenvalue $\beta_{1,1}=\lambda-1$ is simple and changes sign at $\lambda=1$, where all other eigenvalues are at most $\max\{-3,\,1-d\}<0$.
Thus the PES condition holds with $\lambda_0 = 1$ and $c=1$.

\subsubsection*{Self interaction and cross interaction numbers} The self-interaction number $\mathcal{S}_{N}$ given by \autoref{S-definition} equals to
\begin{equation}\label{example-SN-last}
    \mathcal{S}_{N}= \frac{2}{\pi} B\!\left(\frac{N+2}{2},\frac{1}{2}\right) a_{0,N,2}.
\end{equation}
\begin{remark}
If $g_2(0,v)=0$ then $\mathcal{S}_{N}=0 \quad \forall\, N \geq 2$.
\end{remark}

For the cross-interaction number $\mathcal{C}_3$, we use \autoref{C-definition} to get
\begin{equation}\label{example-C2-last}
    \mathcal{C}_3= - \frac{4}{\pi^2} \sum_{m=2} \left[\mathcal{J}_{m,1}^{1,2}\right]^2 \frac{2}{\lambda - m^2} a_{0,2,2}^{2} - \frac{4}{\pi^2} \sum_{m=1} \left[\mathcal{J}_{m,1}^{1,2}\right]^2 \frac{1}{\lambda - dm^2} a_{0,2,1} a_{1,1,2}.
\end{equation}

Since $\left[\mathcal{J}_{m,1}^{1,2}\right]^2$ is always positive, its explicit value is not needed to determine the sign of $\mathcal{C}_3$.
Combining \eqref{example-SN-last} with \eqref{example-C2-last} allows us to analyze the transition type rigorously.

\subsubsection*{Classification of transitions} We will only classify the transitions depending on the coefficients of the Taylor expansion of the nonlinear terms up to order 3.
The classifications depending on the higher order terms can be carried out according to \autoref{main-theorem}.

We note that $\mathcal{S}_2 \ne 0$ if and only if $a_{0,2,2} \neq 0$.
In this case, by \autoref{main-theorem}, the first transition number is
\[
\mathcal{P}= \mathcal{P}_2 = \mathcal{S}_2 = \frac{8}{3\pi} a_{0,2,2},
\]
the critical index is $\mathfrak{m}=2$ and the transition is Type-III.

Next we consider the case
\begin{equation}\label{case-2}
    a_{0,2,2}=0 \text{ } \text{and} \text{ } (a_{0,2,1}a_{1,1,2} \neq 0 \text{ } \text{or} \text{ } a_{0,3,2} \neq 0).
\end{equation}
Since $a_{0,2,2}=0$, we observe that $\mathcal{S}_2 = 0$ and hence $\mathcal{P}_2 = 0$ by \eqref{PN}.

By \eqref{example-SN-last}, the self-interaction number $\mathcal{S}_3$ can be expressed as
\begin{equation*}
    \mathcal{S}_{3}= \frac{3}{4} a_{0,3,2}.
\end{equation*}

Also by \eqref{example-C2-last}, the cross-interaction number $\mathcal{C}_3$ can be expressed as
\begin{equation*}
    \mathcal{C}_3= k(d) a_{0,2,1} a_{1,1,2},
\end{equation*}
where
\begin{equation}\label{kd}
    k(d)= - \dfrac{4}{\pi^2} \sum_{m=1} \left[\mathcal{J}_{m,1}^{1,2}\right]^2 \frac{1}{\lambda - dm^2}.
\end{equation}
By \eqref{J-simplified}, $\mathcal{J}_{m,1}^{1,2}=O(m^{-3})$, so the summands are $O(m^{-8})$ and the series converges absolutely, uniformly on $d \ge 1+\varepsilon$ for each $\varepsilon>0$.
Moreover, at $\lambda=\lambda_0=1$ every denominator satisfies $1-dm^{2}<0$ for $d>1$, so each summand is negative and the leading factor $-4/\pi^{2}$ makes $k(d)>0$, as claimed.

By \eqref{PN}, we have
\begin{equation*}
    \mathcal{P}_3 = \mathcal{C}_3 + \mathcal{S}_3 = k(d) a_{0,2,1} a_{1,1,2} + \frac{3}{4} a_{0,3,2}
\end{equation*}
By \autoref{main-theorem} the conclusions of \autoref{example-main-theorem} hold true.

\subsubsection*{Evaluation of $\mathcal{J}_{m,1}^{1,2}$}

Since $\mathcal{J}_{m,1}^{1,2}$ is required in \eqref{kd}, we provide its explicit calculation below.

In our setting $\gamma = m + 2 (1-j)$, giving three cases:
\begin{enumerate}[label=$\alph*)$]
    \item $j = 0 \implies \gamma_0 = m + 2$.
    \item $j = 1 \implies \gamma_1 = m$.
    \item $j = 2 \implies \gamma_2 = m - 2$.
\end{enumerate}

For $j=0$ we obtain the followings:
\begin{align*}
    &\indic_{\{0\}}(\gamma_j) = 0,  \\
    &\indic_{\{2m\}}(\gamma_j) =
    \begin{cases}
        1, &m=2\\
        0, &\text{otherwise},
    \end{cases} \\
    &\indic_{\mathbb{Z} \setminus \{0, 2m\}}(\gamma_j) =
    \begin{cases}
        0, &m=2\\
        1, &\text{otherwise}.
    \end{cases}
\end{align*}
The cases $j=1$ and $j=2$ are treated analogously.
After substituting these into \eqref{Int-2}, we have
\begin{align*}
    \mathcal{J}_{m,1}^{1,2}
        &= \frac{1}{(2 \ci)^3} \left[0 - \indic_{\{2\}}(m) \pi + \indic_{\mathbb{Z} \setminus \{2\}}(m) \frac{2m\bigl(1-(-1)^m\bigr)}{\ci (m+2) (2-m)} \right]  \\
        &\quad + \frac{1}{(2 \ci)^3} (-2) \left[0 - 0 + 2m \frac{(1-(-1)^m)}{\ci m (-m)} \right] \\
        &\quad + \frac{1}{(2 \ci)^3} \left[ \indic_{\{0\}}(m-2) \pi - 0 + \indic_{\mathbb{Z} \setminus \{0\}}(m-2) \frac{2m\bigl(1-(-1)^m\bigr)}{\ci (m-2) (-m-2)} \right].
\end{align*}

Equivalently,
\begin{equation*}
    \begin{aligned}
        \mathcal{J}_{m,1}^{1,2}= \frac{1}{(2\ci)^3} \left[ (1 - \delta_{m,2} ) \frac{4m\bigl(1-(-1)^m\bigr)}{\ci (m+2) (2-m)}
        + \frac{4m(1-(-1)^m)}{\ci m^2} \right].
    \end{aligned}
\end{equation*}

Finally, $\mathcal{J}_{m,1}^{1,2}$ reduces to
\begin{equation*}
    \mathcal{J}_{m,1}^{1,2}= \frac{m}{2} \left(1 -(-1)^m \right) \left[ (1 - \delta_{m,2} ) \frac{1}{4-m^2} + \frac{1}{m^2} \right].
\end{equation*}
Simplification gives \eqref{J-simplified}.

\subsubsection*{Numerical verification}

We now check the classification of \autoref{example-main-theorem} against numerical solutions of the full system.
The system \eqref{ex-main-eqn}--\eqref{ex-boundary-conditions} is discretized by a sine-spectral method with $200$ modes.
We take $d=2$, for which $k(d)=0.7222$, and measure the solution $\vec{u}=(u,v)^T$ by its projection onto the critical mode $\vec{e}_{1,1}=(0,1)^T\sin x$, which is the variable $u_{c,1}$ of the reduced equation \eqref{reduced-eqn} with $c=1$,
\[
    u_{1,1} = \frac{\langle \vec{u}, \vec{e}_{1,1} \rangle}{\langle \vec{e}_{1,1}, \vec{e}_{1,1} \rangle} = \frac{2}{\pi}\int_0^{\pi} v(x)\,\sin x \, dx .
\]
We consider three scenarios, one for each transition type of \autoref{example-main-theorem}.
By \autoref{example-main-theorem} the transition type depends only on the coefficients listed; the others are arbitrary, and are set to zero in the computations.
\begin{description}
    \item[Scenario A] $a_{0,2,2}=0$, $a_{0,2,1}=1$, $a_{1,1,2}=-2$, $a_{0,3,2}=-3$. Then $\mathfrak{m}=3$ and $\mathcal{P}=-3.694$: Type-I (continuous).
    \item[Scenario B] $a_{0,2,2}=0$, $a_{0,2,1}=2$, $a_{1,1,2}=3$, $a_{0,3,2}=4$. Then $\mathfrak{m}=3$ and $\mathcal{P}=7.333$: Type-II (jump).
    \item[Scenario C] $a_{0,2,2}=1$. Then $\mathfrak{m}=2$ and $\mathcal{P}=8/(3\pi)$: Type-III (mixed).
\end{description}
Scenarios A and B are the two parameter sets of \autoref{table:combined}.

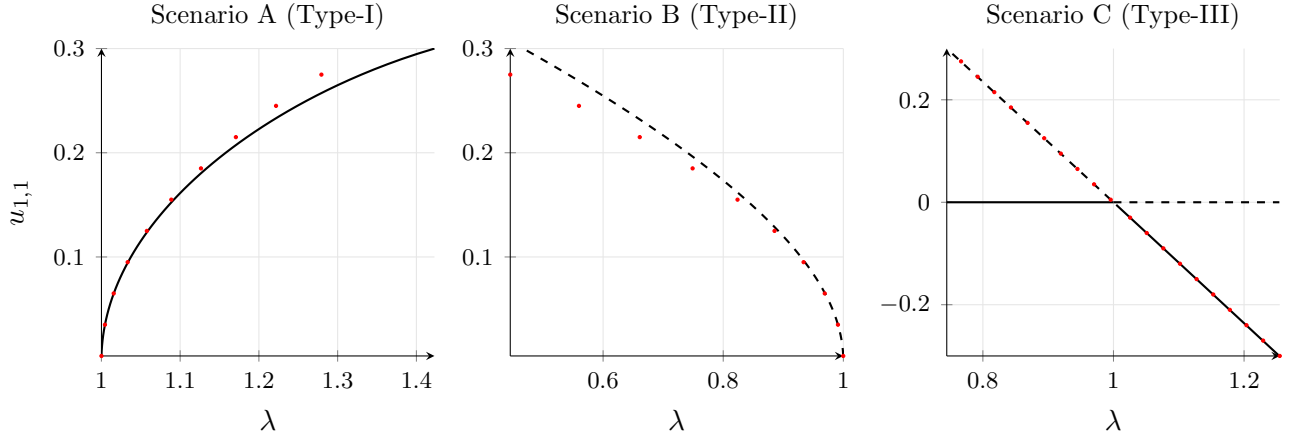
\begin{figure}[H]
    \centering
    \begin{tikzpicture}
        \begin{axis}[width=0.26\textwidth, height=0.24\textwidth, scale only axis, tick label style={font=\scriptsize}, label style={font=\small},
            xlabel={$\lambda$}, ylabel={$u_{1,1}$}, title={\footnotesize Scenario A (Type-I)},
            axis lines=left, grid=major, grid style={gray!20}, title style={yshift=-3pt}]
            \addplot[thick, black] table[x=lam, y=a] {\branchIstable};
            \addplot[only marks, mark=*, mark size=0.6pt, red, each nth point=6]
                table[x=lam, y=a] {\branchIpred};
        \end{axis}
    \end{tikzpicture}\hfill
    \begin{tikzpicture}
        \begin{axis}[width=0.26\textwidth, height=0.24\textwidth, scale only axis, tick label style={font=\scriptsize}, label style={font=\small},
            xlabel={$\lambda$}, title={\footnotesize Scenario B (Type-II)},
            axis lines=left, grid=major, grid style={gray!20}, title style={yshift=-3pt}]
            \addplot[thick, black, dashed] table[x=lam, y=a] {\branchIIunstable};
            \addplot[only marks, mark=*, mark size=0.6pt, red, each nth point=6]
                table[x=lam, y=a] {\branchIIpred};
        \end{axis}
    \end{tikzpicture}\hfill
    \begin{tikzpicture}
        \begin{axis}[width=0.26\textwidth, height=0.24\textwidth, scale only axis, tick label style={font=\scriptsize}, label style={font=\small},
            xlabel={$\lambda$}, title={\footnotesize Scenario C (Type-III)},
            axis lines=left, grid=major, grid style={gray!20}, title style={yshift=-3pt}]
            \addplot[thick, black] table[x=lam, y=a] {\branchIIIstable};
            \addplot[thick, black, dashed] table[x=lam, y=a] {\branchIIIunstable};
            \addplot[thick, black] coordinates {(0.74468055,0) (1,0)};
            \addplot[thick, black, dashed] coordinates {(1,0) (1.2540124,0)};
            \addplot[only marks, mark=*, mark size=0.6pt, red, each nth point=6]
                table[x=lam, y=a] {\branchIIIpred};
        \end{axis}
    \end{tikzpicture}
    \caption{Steady branches at $d=2$: numerical (black; solid stable, dashed unstable) and the analytical estimate from the reduced equation (red).}
    \label{fig:branches}
\end{figure}

Nontrivial steady states are computed by Newton's method.
We prescribe the amplitude $u_{1,1}$ and solve for $\lambda$ together with the solution; since $u_{1,1}\neq 0$ is imposed, the trivial solution is not a root of this system.
Stability is determined from the sign of the largest real part of the eigenvalues of the Jacobian.
\autoref{fig:branches} shows the result.
In each scenario the computed branch agrees near $\lambda_0$ with the prediction $\beta_{1,1}=-\mathcal{P}u_{1,1}^{\mathfrak{m}-1}$ of the reduced equation, shown by the red markers.

\begin{figure}[H]
    \centering
    \begin{tikzpicture}
        \begin{axis}[width=0.36\textwidth, height=0.24\textwidth, scale only axis,
            xlabel={$t$}, ylabel={$u_{1,1}(t)$}, title={\footnotesize Scenario B, $\lambda=0.95$},
            xmin=0, xmax=80, ymin=0, ymax=0.20, ytick={0,0.05,0.1,0.15,0.2}, yticklabel style={/pgf/number format/fixed}, restrict y to domain=-0.2:0.5,
            tick label style={font=\scriptsize}, label style={font=\small},
            axis lines=left, grid=major, grid style={gray!20}, title style={yshift=-3pt}]
            \addplot[thick, gray, dashed, domain=0:80, samples=2] {0.0826};
            \addplot[thick, black] table[x=t, y=a] {\timeIIbelow};
            \addplot[thick, black, densely dotted] table[x=t, y=a] {\timeIIabove};
        \end{axis}
    \end{tikzpicture}\hfill
    \begin{tikzpicture}
        \begin{axis}[width=0.36\textwidth, height=0.24\textwidth, scale only axis,
            xlabel={$t$}, title={\footnotesize Scenario C, $\lambda=1.05$},
            xmin=0, xmax=200, ymin=-0.08, ymax=0.15, restrict y to domain=-0.2:0.5,
            tick label style={font=\scriptsize}, label style={font=\small},
            axis lines=left, grid=major, grid style={gray!20}, title style={yshift=-3pt}]
            \addplot[thick, black] table[x=t, y=a] {\timeIIIleft};
            \addplot[thick, black, densely dotted] table[x=t, y=a] {\timeIIIright};
        \end{axis}
    \end{tikzpicture}
    \caption{Projection $u_{1,1}(t)$ of the solution onto the critical mode for various initial conditions.}
    \label{fig:time}
\end{figure}
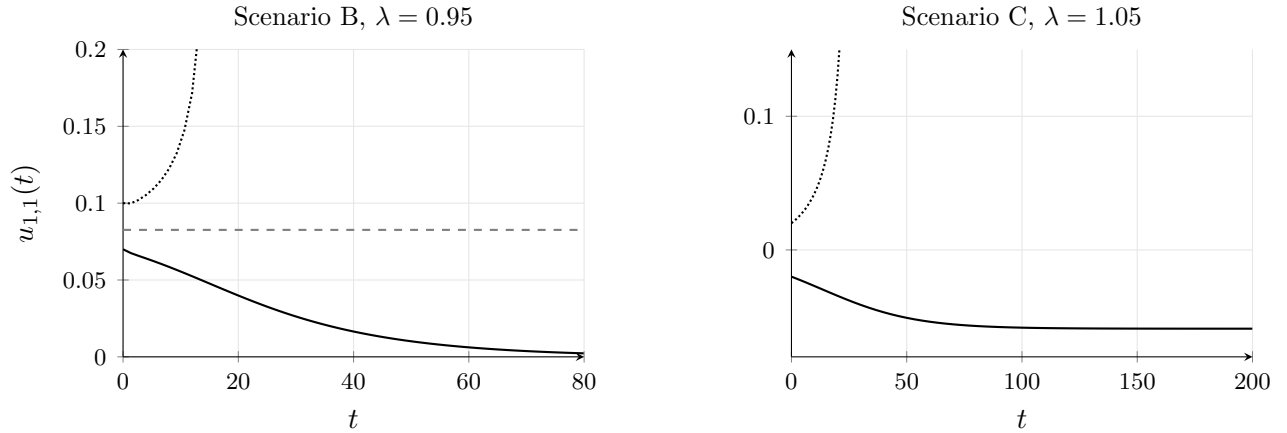

We next integrate the system in time from the initial data $u_{1,1}(0)\sin x$ in the $v$-component, treating the linear part implicitly.
\autoref{fig:time} shows two runs each for Scenarios B and C.
In Scenario B at $\lambda=0.95$, the unstable steady state is at $u_{1,1}=0.0826$ (dashed line).
The solution with $u_{1,1}(0)=0.07$ decays to the trivial state, while the solution with $u_{1,1}(0)=0.10$ leaves the neighborhood, which is the jump.
In Scenario C at $\lambda=1.05$, the solutions with $u_{1,1}(0)=\mp 0.02$ behave in opposite ways.
One settles at $u_{1,1}=-0.0589$, which is the value $-\beta_{1,1}/\mathcal{P}$ on the steady branch, and the other escapes.
This is the splitting of a neighborhood into a continuous region and a jump region which defines a mixed transition.

\subsection{Off-diagonal coupling}\label{subsec:coupled-example}

The reaction matrix \eqref{ex-linear-reaction-matrix} is diagonal, so that the two components of the system are coupled only through the nonlinearity.
We now show that the hypotheses of \autoref{main-theorem} also allow linear cross-coupling, and that the computation remains explicit.
Consider
\begin{equation}\label{coupled-reaction-matrix}
    A=\begin{bmatrix}\lambda & b \\[2pt] b & \lambda\end{bmatrix}, \qquad b \in \mathbb{R},
\end{equation}
in place of \eqref{ex-linear-reaction-matrix}, the rest of the system \eqref{ex-main-eqn}--\eqref{ex-boundary-conditions} being unchanged.
Since $A$ is symmetric, the linear operator is still self-adjoint and the framework of \autoref{sec:formulation-of-the-problem} applies.
The case $b=0$ is the example treated above.

\subsubsection*{Spectrum}

The operator block-diagonalizes over the modes $\sin mx$, the $m$-th block being
\[
    M_m=\begin{bmatrix}\lambda-dm^2 & b \\ b & \lambda-m^2\end{bmatrix} .
\]
Writing
\[
    \theta_m=\frac{(d-1)m^{2}}{2}, \qquad R_m=\sqrt{\theta_m^{2}+b^{2}} ,
\]
the eigenvalues and eigenvectors are
\begin{equation}\label{coupled-spectrum}
    \beta_{m,1;2}=\lambda-\frac{(1+d)m^{2}}{2} \pm R_m, \qquad
    \vec{\alpha}_{m,1} \parallel \begin{bmatrix} b \\ \theta_m+R_m \end{bmatrix}, \quad
    \vec{\alpha}_{m,2} \parallel \begin{bmatrix} b \\ \theta_m-R_m \end{bmatrix},
\end{equation}
the eigenvectors being normalized to unit length.
As $b \to 0$ with $d>1$ these tend to $\vec{\alpha}_{m,1}=(0,1)^{T}$ and $\vec{\alpha}_{m,2}=(1,0)^{T}$, as before.

Since $R_m > |\theta_m|$ for $b \ne 0$, we have $\beta_{m,2}<\beta_{m,1}$ for every $m$.
Also $\beta_{m,1}$ is decreasing in $m$, since the derivative of $-\frac{(1+d)m^{2}}{2}+R_m$ with respect to $m$ is at most $-(1+d)m+|d-1|m<0$.
Hence the largest eigenvalue is $\beta_{1,1}$, it is simple, and the PES condition \eqref{PES-condition} holds with
\begin{equation}\label{coupled-lambda0}
    c=1, \qquad \lambda_0=\frac{1+d}{2}-\sqrt{\frac{(d-1)^{2}}{4}+b^{2}} .
\end{equation}

\subsubsection*{Transition numbers}

Write $\vec{\alpha}=\vec{\alpha}_{1,1}=(\alpha_1,\alpha_2)^{T}$ for the critical eigenvector given by \eqref{coupled-spectrum} at $m=1$.
Since $c=1$ is odd, Remark~\ref{rem:small-cases} gives
\begin{equation}\label{coupled-S2}
    \mathcal{S}_2=\frac{8}{3\pi} \sum_{k+l=2} (\vec{a}_{k,l}\cdot\vec{\alpha})\, \alpha_1^{k} \alpha_2^{l}
    =\frac{8}{3\pi}\Big[ (\vec{a}_{2,0}\cdot\vec{\alpha})\,\alpha_1^{2}
    +(\vec{a}_{1,1}\cdot\vec{\alpha})\,\alpha_1\alpha_2
    +(\vec{a}_{0,2}\cdot\vec{\alpha})\,\alpha_2^{2} \Big],
\end{equation}
and $\mathcal{S}_3=\frac34 \sum_{k+l=3} (\vec{a}_{k,l}\cdot\vec{\alpha})\, \alpha_1^{k} \alpha_2^{l}$.
For the cross term, $\mathcal{J}^{1,2}_{m,1}$ vanishes for even $m$ by \eqref{J-simplified}, so only odd $m$ contribute and
\begin{equation}\label{coupled-C3}
    \mathcal{C}_3=-\frac{4}{\pi^{2}} \sum_{m \text{ odd}} \left[\frac{4}{m(4-m^{2})}\right]^{2}
    \left[ \frac{(1-\delta_{m,1})\,\mathcal{Q}(\vec{\alpha}_{m,1})}{\beta_{m,1}}
    + \frac{\mathcal{Q}(\vec{\alpha}_{m,2})}{\beta_{m,2}} \right],
\end{equation}
where, for a unit vector $\vec{w}$,
\[
    \mathcal{Q}(\vec{w}) = \Big[ \sum_{k+l=2} (\vec{a}_{k,l}\cdot\vec{\alpha}) \big( k\, \alpha_1^{k-1} \alpha_2^{l} w_1 + l\, \alpha_1^{k} \alpha_2^{l-1} w_2 \big) \Big]
    \Big[ \sum_{k'+l'=2} (\vec{a}_{k',l'}\cdot\vec{w})\, \alpha_1^{k'} \alpha_2^{l'} \Big],
\]
all quantities being evaluated at $\lambda=\lambda_0$.
The classification of \autoref{theorem mP} then applies directly: if $\mathcal{S}_2 \ne 0$ the transition is Type-III with $\mathfrak{m}=2$, and if $\mathcal{S}_2=0$ the transition is continuous or jump according to the sign of $\mathcal{P}_3=\mathcal{S}_3+\mathcal{C}_3$.

Unlike the diagonal case, both families of stable modes now contribute to $\mathcal{C}_3$.
At $b=0$ we have $\vec{\alpha}=\vec{\alpha}_{m,1}=(0,1)^{T}$ and $\vec{\alpha}_{m,2}=(1,0)^{T}$, so that
\[
    \mathcal{Q}(\vec{\alpha}_{m,1}) = 2\,a_{0,2,2}^{2}, \qquad \mathcal{Q}(\vec{\alpha}_{m,2}) = a_{0,2,1}\,a_{1,1,2} .
\]
The first vanishes when $\mathcal{S}_2=0$, so the modes $\vec{\alpha}_{m,1}$ drop out and \eqref{coupled-C3} reduces to $k(d)\,a_{0,2,1}a_{1,1,2}$ as in \autoref{example-main-theorem}.

\section{Conclusions}\label{sec:conclusions}
In this work, we have analyzed the first dynamic transition of a general two-component reaction--diffusion system on a bounded interval in one spatial dimension.
This work is part of our ongoing research project on the study of first dynamic transitions of general families of nonlinear partial differential equations.
The emphasis of this research project is on the exact determination of the effects of the nonlinearity of the system on the first dynamic transitions.
\emph{The main goal of this study is the extension of the exact formulation of the first transition in terms of the data of the system, obtained in previous studies \cite{csengul2022dynamic, csengul2023dynamic, csengul2024first} which focused on single-equation cases, to systems of coupled equations}.

During this work we encountered many possible future related directions which we list below.
\begin{enumerate}
    \item The principle of exchange of stabilities condition may be satisfied by more than one eigenvalue. This is crucial for understanding pattern selection in reaction diffusion systems.
    \item The nonlinear terms may depend on the derivatives of the unknown variables as well. For a single equation this has been carried out; for systems it is open.
    \[
    G(u,v) = \sum_{\abs{\alpha_u} + \abs{\alpha_v} = 2}^{\infty} \vec{a}_{\alpha_u, \alpha_v} (D u)^{\alpha_u} (D v)^{\alpha_v},
    \]
    where $\vec{a}_{\alpha}$ are real coefficients and $\alpha_u, \alpha_v \in \mathbb{Z}_{\geq 0}$ are multi-indices. This extends reaction diffusion systems to more general nonlinear partial differential equations.
    \item The removal of the self-adjointness hypothesis, which is what separates the present results from the classical two-component models such as the Brusselator, Gray--Scott, Schnakenberg and Gierer--Meinhardt systems.
    In that case the eigenvectors and the adjoint eigenvectors differ, and the projections in the reduction must be taken with respect to the adjoint eigenvectors, as in the single-equation case \cite{csengul2023dynamic}.
    We plan to carry this out in a forthcoming paper.
    \item We will consider the case where the eigenfunctions consist of pure cosine modes in a future work. This is often the case where the system \eqref{main-eqn} is considered with homogeneous Neumann boundary conditions which have many interesting applications in practice.
\end{enumerate}

\section*{CRediT Authorship Contribution Statement}
\textbf{Taylan Şengül:} Conceptualization, Methodology, Supervision, Writing-review \& editing - \textbf{Burhan Tiryakioglu, Esmanur Yıldız Akıl:} Conceptualization, Formal analysis, Investigation, Writing-original draft

\section*{Conflicts of Interest}
This work does not have any conflicts of interest.

\section*{Use of AI Tools Declaration}
The authors declare they have not used Artificial Intelligence (AI) tools in the creation of this article.

\section*{Acknowledgements}
We thank the two anonymous referees for their careful reading and constructive comments, which led to a clearer and more complete paper.
This study was supported by the Scientific and Technological Research Council of Türkiye (TÜBİTAK), Science Fellowships and Grant Programs Department (BİDEB), under the 2211-A National PhD Scholarship Program.

% \nocite{*}  % uncomment to list uncited bib entries as well
\bibliographystyle{plain}
\bibliography{1d2v}

\end{document}